\newcount\remarkcc

\catcode`@=11

\newdimen\iiiw
\newdimen\iiih
\def\iiib{_{\kern0pt
	\raise1pt\hbox to 0pt{\vrule width\iiiw height\iiih depth0pt
		\dimen0=\iiiw \divide\dimen0 by 2 \kern-\dimen0
		\raise-\dimen0\hbox to 0pt{\vrule width\iiih height\iiiw depth0pt\hss}
		\hss}
	\kern\iiiw
	\kern 2pt
	}}
\def\iiia{_{\kern0pt
	\raise1pt\hbox to 0pt{\vrule width\iiiw height\iiih depth0pt
		\hss}
	\kern\iiiw
	\kern 2pt
	}}

\def\fsize{12pt}

	\documentclass[\fsize,a4paper,twoside]{article}
	\usepackage{amsfonts,amsmath,amssymb,color,epsfig,cite}
\def\endproof{\mbox{}\nobreak\hfill\nobreak$\Box$\par\vspace{2ex}}

\definecolor{ccred}{rgb}{1,0,0}
\definecolor{ccgreen}{rgb}{0,0.6,0}
\definecolor{ccnix}{rgb}{1,1,1}
\definecolor{ccwhite}{rgb}{1,1,1}
\definecolor{ccblue}{rgb}{0.5,0.5,1}
\definecolor{cclila}{rgb}{1,0,1}
\definecolor{ccblack}{rgb}{0,0,0}

\newcommand{\nix}[1]{\textcolor{ccnix}{#1}}
\newcommand{\white}[1]{\textcolor{ccnix}{#1}}

\newdimen\xind
\def\versiony#1#2%
  {\leavevmode\hbox to 0pt{\vbox to 0pt{{\vs#2\tiny ver #1}\vss}\hss}}

\def\bbig#1{{\hbox{$\left#1\vbox to 10.9pt{}\right.\n@space$}}}

\def\rect#1#2#3#4
  {\setbox0=\vbox{\vrule width#3 height#4}%
   \hbox to 0pt{\hs#1\raise#2\box0\hss}%
  }

\def\RaCut{\hbox to 9pt{$\Ra%
   \white{%
     \rect{-9pt}{4.6pt}{3pt}{2pt}%
     \rect{-9pt}{-0.6pt}{3pt}{2pt}%
     \rect{-8pt}{1pt}{2.3pt}{4pt}%
     \rect{-7pt}{2.15pt}{3pt}{1.7pt}%
   }$\hss}}%

\newdimen\hsizee
\def\rucky#1#2#3{{\hsizee=\hsize\advance\hsizee by -#1\par\penalty5000\leavevmode\hangindent=#1\hangafter=1\hbox to #1{#2}#3\par}}

\def\rulenix{\vrule width0pt depth0pt height0pt}
\let\nix=\rulenix

   \let\hs=\hskip
   \let\vs=\vskip

   \def\es{\enspace }
   \def\ess{\es\es}

   \def\1{\ifmmode\mskip1.5mu\else\hskip0.06em\relax\fi} 
   \def\,{\ifmmode\mskip\thinmuskip\else\hskip0.167em\relax\fi}
   \def\>{\ifmmode\mskip\medmuskip\else\hskip0.222em\relax\fi}
   \def\;{\ifmmode\mskip\thickmuskip\else\hskip0.278em\relax\fi}
   \def\!{\ifmmode\mskip-\thinmuskip\else\hskip-0.167em\relax\fi}

   \def\hb{\hs-1.3cm}

  \def\iMax#1#2{\ifnum #1>#2 #1\else #2\fi}
  \def\iWidth#1#2
       {\ifmmode%
          \setbox0=\hbox{$#2$}#1=\wd0%
        \else%
          \setbox0=\hbox{#2}#1=\wd0%
        \fi%
       }
  \def\iHeight#1#2
       {\ifmmode%
          \setbox0=\hbox{$#2$}#1=\ht0%
        \else%
          \setbox0=\hbox{#2}#1=\ht0%
        \fi%
       }
  \def\iDepth#1#2
       {\ifmmode%
          \setbox0=\hbox{$#2$}#1=\dp0%
        \else%
          \setbox0=\hbox{#2}#1=\dp0%
        \fi%
       }
  \def\iFrac#1#2#3{\multiply #1 by #2\divide #1 by #3}
  \def\iKlamm#1#2{\hbox{$#1{\vbox to#2{}}$}}

  \def\iKlammSysA
    {%
      \def\iLeft##1{\left(##1\right.}%
      \def\iMid##1{\left|##1\right.}%
      \def\iRight##1{\left)##1\right.}%
      \def\iSubSkip{-1.5pt}%
    }%
  \def\iKlammSysB
    {%
      \def\iLeft##1{\left\langle##1\right.}%
      \def\iMid##1{,##1}%
      \def\iRight##1{\left\rangle##1\right.}%
      \def\iSubSkip{-1pt}%
   }%

  \def\iKlammSysI%
    {%
      \def\iMid##1{,##1}%
      \def\iSubSkip{-1pt}%
    }%
  \def\iKlammSysIoo
    {\iKlammSysI\def\iLeft##1{\left]\right.}\def\iRight##1{\left[\right.}}
  \def\iKlammSysIoa
    {\iKlammSysI\def\iLeft##1{\left]\right.}\def\iRight##1{\left]\right.}}
  \def\iKlammSysIao
    {\iKlammSysI\def\iLeft##1{\left[\right.}\def\iRight##1{\left[\right.}}
  \def\iKlammSysIaa
    {\iKlammSysI\def\iLeft##1{\left[\right.}\def\iRight##1{\left]\right.}}

  \newdimen\iDa
  \newdimen\iDb
  \newdimen\iDc

  \def\iKlammer#1#2
      {%
        \iWidth\iDa{#1}%
        \iWidth\iDb{#2}%
        \iDa=\iMax\iDa\iDb%
        \iFrac\iDa{1}{80}
        \iDc=\iMax\iDa\p@%
        \iHeight\iDa{#1}%
        \iHeight\iDb{#2}%
        \iDa=\iMax\iDa\iDb%
        \advance\iDa by 1pt
        \iKlamm\iLeft\iDa\hskip\iDc\hskip-1.5pt%
          #1%
        \hskip\iDc\iKlamm\iMid\iDa\hskip\iDc\hskip-1.5pt%
          #2%
        \hskip\iDc\hskip-.1pt\iKlamm\iRight\iDa%
      }

  \def\iSub#1{\hskip\iSubSkip{\vrule width0pt depth3pt}_{#1}}

  \def\iKlammerA{\iKlammSysA\iKlammer}
  \def\iKlammerB{\iKlammSysB\iKlammer}

  \def\iBlock#1
      {%
        \iWidth\iDa{#1}%
        \iFrac\iDa{1}{90}%
        \iDc=\iMax\iDa\p@%
         \hskip\iDc%
          #1%
         \hskip\iDc%
      }
  \def\iBlockR#1
      {%
        \iWidth\iDa{#1}%
        \iWidth\iDb{\,}%
        \iFrac\iDa{1}{40}%
        \iDc=\iMax\iDa\iDb%
          #1%
         \hskip\iDc%
      }

  \newdimen\iUserA
  \newdimen\iUserB

  \def\pboxl#1#2{\hbox to #1{#2\hss}}
  \def\pboxc#1#2{\hbox to #1{\hss#2\hss}}
  \def\pboxr#1#2{\hbox to #1{\hss#2}}

  \def\xhbox#1#2%
        {\ifmmode
           \setbox0=\hbox{$#2$}\ht0=#1\box0%
         \else%
           \setbox0=\hbox{#2}\ht0=#1\box0%
         \fi%
        }
  \def\xhdbox#1#2#3%
        {\ifmmode
           \setbox0=\hbox{$#3$}\ht0=#1\dp0=#2\box0%
         \else%
           \setbox0=\hbox{#3}\ht0=#1\dp0=#2\box0%
         \fi%
        }

  \def\bNoHtDp#1{{\def\xa{}\ifmmode\def\xa{$}\fi\setbox0=\hbox{\xa #1\xa}\ht0=0pt\dp0=0pt\box0}}%

  \def\os#1#2{\ifmmode\setbox#1=\hbox{$#2$}\else\setbox#1=\hbox{#2}\fi}
  \def\osD#1#2{\os0{#1}\os1{#2}\dp0=\dp1\box0}
  
  \def\osH#1#2{\os0{#1}\os1{#2}\ht0=\ht1\box0}

\def\dwh#1#2{{
  \setbox0=\hbox{$#2$}
  \setbox1=\hbox{\vrule width\wd0}
  \setbox2=\hbox{\hs#1$\wh{\box1}$}
  \iDa=\ht0\advance\iDa by -6pt
  \wd2=0pt
  \raise\iDa\copy2
  \advance\iDa by 1.5pt
  \raise\iDa\copy2
  \box0
}}

  \newdimen\ullh

   \def\f#1#2{{#1 \over #2}}
   \def\ul{\underline}
   \def\ull#1{\setbox0=\hbox{$#1$}\ullh=\wd0\advance\ullh by -3pt\wd0=\ullh\ul{\box0}\kern3pt}
   \def\xwt#1{\iHeight{\dimen0}{#1}\xhbox{\dimen0}{\wt{#1}}}
   \def\Int{\mathop{\rm Int}}
   \def\sm{\setminus}

   \def\wcl#1
      {\overline{#1}^{\hs0.8pt w}}

   \def\joinrell{\mathrel{\mkern-6mu}}%

   \def\emp{\emptyset}

   \def\Lral{\Longleftrightarrow}

   \def\nRal{$\hbox to 0pt{\kern3pt$\not$\hss}$\Ral}

   \def\Ral{\Longrightarrow}
   \def\Ra{\Rightarrow}

   \def\dLral{\ifmmode \ess\Lral\ess \else $\ess\Lral\ess$ \fi}

   \def\Ge{\varepsilon}

\begin{document} 

\newcounter{tcnt}[section]
\renewcommand{\thetcnt}{\thesection.\arabic{tcnt}}

\newtheorem{remark}[tcnt]{Remark}
\newtheorem{example}[tcnt]{Example}
\newtheorem{assumption}[tcnt]{Assumption}
\newtheorem{convention}[tcnt]{Convention}
\newtheorem{problem}[tcnt]{Problem}

\newtheorem{definition}[tcnt]{Definition}
\newtheorem{lemma}[tcnt]{Lemma}
\newtheorem{corollary}[tcnt]{Corollary}
\newtheorem{proposition}[tcnt]{Proposition}
\newtheorem{theorem}[tcnt]{Theorem}

\def\pboxl#1#2{\hbox to #1{#2\hss}}
\def\pboxc#1#2{\hbox to #1{\hss#2\hss}}
\def\pboxr#1#2{\hbox to #1{\hss#2}}

\def\xhh#1#2{\iHeight{\dimen0}{#2}\xhbox{\dimen0}{#1}}
\def\xwt#1{\xhh{\widetilde  #1}{#1}}

\author{D.~Wegner\thanks{Department of Mathematics, Humboldt-Universit\"at zu Berlin, E-Mail: dwegner@math.hu-berlin.de}}
\title{Existence for coupled pseu\discretionary{-}{}{}domono\discretionary{-}{}{}tone--strongly monotone systems and application to a Cahn--Hilliard model with elasticity}
\date{}

\maketitle

\begin{abstract}
	A system of two operator equations is considered --
		one of pseu\discretionary{-}{}{}domono\discretionary{-}{}{}tone type and the other of strongly monotone type --
		both being strongly coupled.
	Conditions are given that allow to reduce the solvability of this system to a
		single operator equation for a pseu\discretionary{-}{}{}domono\discretionary{-}{}{}tone mapping.
	This result is applied to a coupled system consisting of
		a parabolic equation of forth order in space of Cahn--Hilliard type and
		a nonlinear elliptic equation of second order 
		to a quasi-steady mechanical equilibrium.
	Using an appropriate notation of weak solutions and a framework for evolution equations
		developed by Gr\"oger~\cite{konni},
		the system is reduced to a single parabolic operator equation and the existence of solutions
		are shown under restrictions on the strength of the coupling.
%
%
\end{abstract}

\section{Introduction}

In this paper the pseudomonotonicity of special compositions	of nonlinear operators is shown.
More specifically, we consider the system
	\begin{eqnarray*} 
		A (»x,»y)	& = &	x_0^*,	\\
		B (»x,»y)	& = &	y_0^*
	\end{eqnarray*}
	for operators $ A : X \times Y  \rightarrow   X ^*$ and $ B : X \times Y  \rightarrow   Y ^*$ on reflexive Banach spaces $ X $ and $ Y $.
Assume that for every $»x\in X $ the mapping $»Bx:= B (»x,\>.\> ): Y \rightarrow   Y ^*$ is uniquely invertible
	and define $ R »x:=»Bx^{-1}( y_0^*)$.
Then the given system is solvable if and only if $ A (»x, R »x)= x_0^*$ admits a solution.\\
We provide sufficient conditions that ensure the pseudomonotonicity of the mapping $ S »x:= A (»x, R »x)$.
Then existence results for this system can be obtained from classical theory of pseu\discretionary{-}{}{}domono\discretionary{-}{}{}tone operators.
To this end, we introduce the subclass of semimonotone operators
	(which is a variant of a respective subclass of pseu\discretionary{-}{}{}domono\discretionary{-}{}{}tone operators considered in
	\cite{lions,zeid2b,papa,deim}).
The operators of this subclass enjoy a mixture of monotonicity and of compactness properties. 
This can be seen as a generalization of those differential operators that
	are monotone in the highest order terms and	compact in the terms of lower order.\\
The conditions we presume in order to prove the pseudomonotonicity of $ S $
	consist of the strong monotonicity of $ B $ in $»y$,
	the semimonotonicity of $ A $ in $»x$,
	and further assumptions on the coupling of both equations of the system.
The latter include respective Lipschitz conditions.
Furthermore, we require that when splitting $ A $ into a monotone and a compact part
	the composition operator $ S $ still inherits the monotonicity property of $ A $.
This leads to a restriction on the influence both parts of the system may exert on each other 
	and is given as a condition on the Lipschitz and strong monotonicity constants.

{\def\xa#1#2#3{&\pboxc{0cm}{\pboxc{0.87\hsize}{$#1$}\pboxl{0.04\hsize}{#2}\pboxl{0.09\hsize}{$#3$}}&}
In our application this reduction technique is applied to
	a model of phase-separation in a binary mixture incorporating elastic effects.
To be more specific, we consider
	on a time interval ${\cal T}$ 
	and on a domain $\Omega $, with ${\Gamma_D}$ and ${\Gamma_N}$ being disjoint parts of the boundary,
	the following parabolic equation of forth order in space of Cahn--Hilliard type
	coupled to an elliptic equation accounting for elastic effect given by
	\begin{eqnarray*} 
		\xa{	\partial_t »u - \mathop{\rm div}(»M\nabla(\mu \partial_t »u+»w)) \>=\> 	0					}{on}{{\cal T}\times\Omega ,}	\\
		\xa{	»w  =  \varphi '(»u) - \mathop{\rm div}( b_1 (»u,\nabla»u, e ))+ b_2 (»u,\nabla»u, e )							}{on}{{\cal T}\times\Omega ,}	\\
		\xa{	\mathop{\rm div} b_0 (»u,\nabla»u, e ) \>=\>  0,\enspace \enspace \enspace 		 e =\epsilon({\bf u}):=\frac12(D{\bf u}+D{\bf u}^t)	}{on}{{\cal T}\times\Omega ,}
	\end{eqnarray*}
	together with the boundary and initial conditions
	\begin{eqnarray*} 
		\xa{	»M\nabla(\mu \partial_t »u+»w)\cdot\vec n	\>=\> 0,\enspace \enspace \enspace    b_1 (»u,\nabla»u, e )\cdot\vec n\>=\>  0	}{on}{{\cal T}\times\partial\Omega  ,}	\\
		\xa{	b_0 (»u,\nabla»u, e )\vec n 			\>=\>   0				}{on}{{\cal T}\times{\Gamma_N},}	\\
		\xa{	{\bf u}  \>=\>   0											}{on}{{\cal T}\times{\Gamma_D},}	\\
		\xa{	»u(0)  \>=\>   u_0 										}{on}{\Omega .}
	\end{eqnarray*}
These equations model
	the mass balance for the concentration $»u$ of one of the components,
	the related chemical potential $»w$,
	and a quasi-steady mechanical equilibrium, respectively,
	with $ b_0 $ being the stress tensor
	which depends in a nonlinear way on $»u,\,\nabla»u$ and on the linearized strain tensor $ e $.
The latter is given as the symmetric part of the derivative of the displacement ${\bf u}$.
Furthermore, $»M$ is the (constant) mobility matrix,
	the functions $ b_1 $ and $ b_2 $ together with convex functional $\varphi $
	determine the chemical potential $»w$ and thus model the behavior of the material.
Note that $ b_0 , b_1 $ and $ b_2 $ may explicitly depend on $(t,x)\in{\cal T}\times\Omega $,
	which is suppressed in the notation to enhance the readability. 
The constant $\mu $ is non-negative.
If it is strictly positive, then the model includes additional contributions to the diffusion flux
	resulting from the concept of microforces, cf.~Fried, Gurtin~\cite{fried1,fried2}
	and Gurtin~\cite{gurtin}.

We prove the existence of solutions in an appropriate weak sense.
For this purpose, we make use of a general framework for evolution equations by Gr\"oger~\cite{konni},
	which allows to include suitable (possibly degenerate) linear operators inside the time derivative.
Then, using our general result the coupled system can be reduced to a single parabolic operator equation
	involving a pseu\discretionary{-}{}{}domono\discretionary{-}{}{}tone operator.
To this end, 
	we have to ensure the aforementioned assumptions on the coupling.
This is done with the help of result on $W^{1,p}$ regularity for some $p>2$ for the solution ${\bf u}$
	to the mechanical equilibrium.

For different models of Cahn--Hilliard type for phase separation coupled to elastic effects
	and related existence existence result exemplarily we refer to \cite{miranville,garcke1,sprekels}.
In~\cite{pawlow} the elastic effects are not assumed to be quasi-steady.
This leads to a coupled system of parabolic-hyperbolic type.
A model which incorporates a damage process was considered in~\cite{heinemann}.

The remainder of this paper is organized as follows:
In Section~\ref{s2:} we introduce our notion of semimonotone operators and show that they form
	a subclass of all pseu\discretionary{-}{}{}domono\discretionary{-}{}{}tone mappings.
Further, we provide conditions on $ A $ and $ B $ such that $ S $ is semimonotone.
Section~\ref{s3:} gives a short introduction into the approach to evolution equations developed by
	Gr\"oger~\cite{konni} and states a corresponding existence result.
In Section~\ref{s4:} the results of the preceding sections are applied to the model above of
	phase-separation with elastic effects.
We introduce an appropriate notion of weak solutions and give conditions on the functions
	$ b_0 , b_1 , b_2 $ and $\varphi $ that are used in order to prove the existence of solutions.

\section{Semimonotone operators}
\label{s2:}

This section introduces the class semimonotone operators and shows that it is a subset of all pseu\discretionary{-}{}{}domono\discretionary{-}{}{}tone operators.
Conditions are given that ensure the composition of two operators with special properties to be semimonotone.
In Section~\ref{s4:} we use this result to reformulate an elliptic-parabolic system as a single evolution
	equation of pseu\discretionary{-}{}{}domono\discretionary{-}{}{}tone type.
For this equation we derive an existence result from the classical theory of pseu\discretionary{-}{}{}domono\discretionary{-}{}{}tone operators.

Before starting with our analysis, let us fix some notations.
For a Banach space $ X $, we denote by $ ||\iBlock{.}||  _{ X } $ its norm, its dual space by $ X ^*$ and
	with $ \iKlammerB{.}{.}  \iSub{ X }: X ^*\times X \rightarrow  {{\mathbb R}}$ its dual pairing.
In this paper, we will only consider real Banach space.
$ X _\omega $ indicates the spaces $ X $ equipped with its weak topology.
If it is clear from the context, we simply write $ ||\iBlock{.}|| $ and $ \iKlammerB{.}{.} $
	for $ ||\iBlock{.}||  _{ X } $ and $ \iKlammerB{.}{.}  \iSub{ X }$, respectively.
Moreover, the (in general multi-valued) duality mapping of $ X $ is given by $ J_{ X } \subset X \times X ^*$.
Here and below, we identify mappings with their graphs and, occasionally,
	singletons $\{x\}$ with $x$ itself.
For a Hilbert space $ H $ we denote by $ \iKlammerA{.}{.}  \iSub{ H }$ its inner product.
Then $ J_{ H } $ coincides with the canonical isomorphism from $ H $ onto $ H ^*$
	given by Riesz's theorem.
The identity mapping of set $M$ regarded as an operator from $M$ into some superset $M'\supset M$
	is written as ${\rm Id}_{M\rightarrow   M'}$.
Finally, for sets $ M_1 , M_2 , M_3 $, $x\in M_1 $ and $F: M_1 \times M_2 \rightarrow   M_3 $
	we write $F_x: M_2 \rightarrow   M_3 $ for the mapping $ M_2 \ni y\mapsto  F(x,y)$.

Now, let $ X $ and $ Y $ be real, reflexive Banach spaces.
We start by recalling the definition of pseu\discretionary{-}{}{}domono\discretionary{-}{}{}tone operators.

\begin{definition}[{\rm $ T $--pas, pseu\discretionary{-}{}{}domono\discretionary{-}{}{}tone operators}]\label{..}
	Let $ T : X \rightarrow   X ^*$ be an operator.
	A sequence $( x_n )_{n\in{{\mathbb N}}}$ in $ X $ will be called a $ T $--pas if $( x_n )$
	converges weakly in $ X $ to an element $»x\in X $ and it holds
	\[  \mathop{\overline {\lim}}_{n\mathop{\rightarrow }\infty}  \iKlammerB{ T x_n }{ x_n -»x}  \>\le\>  0.  \]
	Furthermore, $ T $ is said to be pseu\discretionary{-}{}{}domono\discretionary{-}{}{}tone if 
	for every $ T $--pas $( x_n )_{n\in{{\mathbb N}}}$ converging weakly to $»x\in X $,
	\[   \iKlammerB{ T »x}{»x-v}  \>\le\>   \mathop{\underline {\lim}}_{n\mathop{\rightarrow }\infty} \iKlammerB{ T x_n }{ x_n -v} \]
	holds for every $v\in X $.
\end{definition}

Our notational shortcut of a $ T $--pas stands for a 'pseudomonotonously active sequence'.
The definition of pseudomonotonicity follows Zeidler~\cite{zeid2b}.
Note that the original definition of Br\`ezis involves nets instead of sequences
	and requires the operator to satisfy a certain boundedness condition.

\begin{remark}\label{r1:Doni}
	If $ T : X \rightarrow   X ^*$ is pseu\discretionary{-}{}{}domono\discretionary{-}{}{}tone and if $( x_n )$ is $ T $--pas with limit $»x$,
		then by choosing $»v=»x$ we obtain $\mathop{\underline {\lim}}_n \iKlammerB{ T x_n }{ x_n -»x} \ge0$
		and hence $\lim_n \iKlammerB{ T x_n }{ x_n -»x} =0$.
\end{remark}

\begin{definition}\label{}
	For arbitrary vector spaces $U$ and $V$, $ u_0 \in U$ and for any multi-valued operator $T\subset U\times V$
		the translation $\mathfrak{T}_{ u_0 }T\subset U\times V$ of \,$T$ is given by
		$(\mathfrak{T}_{ u_0 }T)u:=\mathfrak{T}_{ u_0 }Tu:=T(u- u_0 )$ for $u\in U$.
\end{definition}

An important consequence of the class of pseu\discretionary{-}{}{}domono\discretionary{-}{}{}tone operators is its closedness under summation and translation.

\begin{proposition}\label{p2:PsmSum}
	If $ x_0 \in X $ and if the operators $ T , T_1 , T_2 : X \rightarrow   X ^*$ are pseu\discretionary{-}{}{}domono\discretionary{-}{}{}tone, so are $ T_1 + T_2 $ and $\mathfrak{T}_{ x_0 } T $.
\end{proposition}

{\em Proof}. 
	For the pseudomonotonicity of $ T_1 + T_2 $ see~\cite{zeid2b}, Prop~27.6, p.~586.
	Let $( x_n )$ a $\mathfrak{T}_{ x_0 } T $--pas with weak limit $»x$ and $»v$ be arbitrary.
	Then $ y_n := x_n - x_0 $ is a $ T $--pas with limit $»y:=»x- x_0 $ and by the pseudomonotonicity of $ T_1 $ we get
		for $»u:=»v- x_0 $ that
	\[
		 \iKlammerB{\mathfrak{T}_{ x_0 } T »x}{»x-»v} 
			\>=\> 		 \iKlammerB{ T »y}{»y-»u} 
			\>\le\> 		\mathop{\underline {\lim}}_{n\mathop{\rightarrow }\infty}   \iKlammerB{ T y_n }{ y_n -»u} 
			\>=\> 		\mathop{\underline {\lim}}_{n\mathop{\rightarrow }\infty}   \iKlammerB{\mathfrak{T}_{ x_0 } T x_n }{ x_n -»v} 
	\]
	which finishes the proof.
\endproof

\begin{definition}\label{}
	Let $ L :D( L )\rightarrow   X ^*$ be a linear, closed operator with domain $D( X )$ dense in $ X $.
	We set $»Z:=D( L )$ and equip its with the graph norm of $ L $, i.e.
	\[
			 ||\iBlock{»x}||  _{Z} 	:=	( ||\iBlock{»x}||  _{ X } ^2 +  ||\iBlock{ L »x}||  _{ X ^*} ^2)^{1/2}.
	\]
	An operator $ T : X \rightarrow   X ^*$ is said to be pseu\discretionary{-}{}{}domono\discretionary{-}{}{}tone with respect to $ L $,
		if $I^* T  I:Z\rightarrow   Z^*$ is pseu\discretionary{-}{}{}domono\discretionary{-}{}{}tone, whereas $I:={\rm Id}_{Z\rightarrow   X }$ is the identity regarded
		as a mapping from $Z$ into $ X $.\\
	$ T : X \rightarrow   X ^*$ is called radially continuous in $»x\in X $ if the mapping
		$t\mapsto  \iKlammerB{ T (»x+t»v)}{»v} $ from ${{\mathbb R}}$ into itself is continuous in $t=0$ for every $»v\in X $.
	Finally, we call $ T : X \rightarrow   X ^*$ coercive with respect to $ x_0 \in X $ if
	\[
		\lim_{ ||\iBlock{»x}|| \mathop{\rightarrow }\infty}  \iKlammerB{ T »x}{»x- x_0 }  \>=\>  0.
	\]
\end{definition}

Pseudomonotone operators occurring in PDEs often have a special structure:
	a monotone part (usually terms of highest order) together with a
	compact perturbation (lower order terms).
The following notion generalizes this behavior.

\begin{definition}[{\rm Semimonotone operators}]\label{d1:Semi} \def\xa#1#2#3{#1&\hskip3mm& \hbox to 12.5cm{$\displaystyle #2$\hfil$#3$}}
	We call an operator $ T : X \rightarrow   X ^*$ semimonotone if $ T $ has the form $ T »x=\xwt{ T }(»x,»x)$
	for a mapping $\xwt{ T }: X \rightarrow   X ^*$ satisfying the conditions:
	\begin{eqnarray*} 
	\xa{(S1)}{   \iKlammerB{\xwt{ T }(»x,»x)-\xwt{ T }(»y,»x)}{»x-»y}  \>\ge\>  0                                            }{  \forall»x,»y\in X },\\
	\xa{(S2)}{  y_n \mathop{\rightharpoonup}»y \mbox{ \>is a $ T $--pas} \enspace \enspace \Longrightarrow \enspace \enspace   \xwt{ T }(»x, y_n )\mathop{\rightharpoonup}\xwt{ T }(»x,»y)		}{  \forall»x\in X  },\\
	\xa{(S3)}{  y_n \mathop{\rightharpoonup}»y \mbox{ \>is a $ T $--pas} \enspace \enspace \Longrightarrow \enspace \enspace    \iKlammerB{\xwt{ T }(»x, y_n )}{ y_n -»y}  \mathop{\rightarrow }0}{  \forall»x\in X  },\\
	\xa{(S4)}{  »x\mapsto \xwt{ T }(»x,»y) \mbox{ \> is radially continuous in the point $»x=»y$}                         }{  \forall»y\in X  }.
	\end{eqnarray*}
	In this case, $\xwt{ T }$ is called a semimonotone representative of $ T $.
\end{definition}

\begin{remark}\label{..}
	Different authors denote different classes of operators as being semimonotone.
	Deimling~\cite{deim}, Zeidler~\cite{zeid2b} and Hu/Papageorgiou~\cite{papa}
		use definitions which are more restrictive than~\ref{d1:Semi}
		as well as Lions~\cite{lions} and his operators of 'variational type'.
	We use Definition \ref{d1:Semi} instead, since it is more simple and more general,
		but nevertheless collects all the properties we need.
\end{remark}

The following proposition shows that semimonotone operators are indeed pseu\discretionary{-}{}{}domono\discretionary{-}{}{}tone.

\begin{proposition}\label{p1:SemiPsm}
	If $ T : X \rightarrow   X ^*$ is a semimonotone operator, then $ T $ is pseu\discretionary{-}{}{}domono\discretionary{-}{}{}tone.
\end{proposition}

{\em Proof}. 
	Let $( x_n )$ be a $ T $--pas with $ x_n \mathop{\rightharpoonup}»x$, $»v\in X$ and $\xwt{ T }$ be a semimonotone representative of $ T $.
	We put $ w_t :=»x+t(v-»x)$ for $0<t\le1$.
	The monotonicity condition (S1) applied to $ x_n $ and $ w_t $ implies
	\[
		 \iKlammerB{\xwt{ T }( x_n , x_n )-\xwt{ T }( w_t , x_n )}{ x_n -»x+»x- w_t }  \>\ge\>   0.
	\]
	With $»x- w_t =t(»x-»v)$, this can be rewritten as
	\[
		t\,  \iKlammerB{ T x_n }{»x-v}  \enspace \ge\enspace   -\> \iKlammerB{ T x_n }{ x_n -»x}  »+  \iKlammerB{\xwt{ T }( w_t , x_n )}{ x_n -»x}  »+ t\,  \iKlammerB{\xwt{ T }( w_t , x_n )}{»x-v} .
	\]
	Passing to the limit inferior on both sides and using (S2), (S3) and the fact that $( x_n )$
		is a $ T $--pas, we end up with
	\[
		t \!\mathop{\underline {\lim}}_{n\mathop{\rightarrow }\infty}  \iKlammerB{ T x_n }{»x-v}   \>\ge\>    t\, \iKlammerB{\xwt{ T }( w_t ,»x)}{»x-v} .
	\]
	Now we divide by $t$ and pass with $t\mathop{\rightarrow }0$ to the limit in order to obtain
	$
		\mathop{\underline {\lim}}_{n}  \iKlammerB{ T x_n }{»x-v}   \>\ge\>     \iKlammerB{ T »x}{»x-v}   
	$
	by the radial continuity (S4). This inequality together with $\lim_{n}  \iKlammerB{ T x_n }{ x_n -»x} =0$
		(cf.\ Remark~\ref{r1:Doni}) yields
	\[
		\mathop{\underline {\lim}}_{n\mathop{\rightarrow }\infty}  \iKlammerB{ T x_n }{ x_n -v} 
				\enspace \ge\enspace    	\mathop{\underline {\lim}}_{n\mathop{\rightarrow }\infty}  \iKlammerB{ T x_n }{ x_n -»x}  »+ \mathop{\underline {\lim}}_{n\mathop{\rightarrow }\infty}  \iKlammerB{ T x_n }{»x-v}   
				\enspace \ge\enspace  	 \iKlammerB{ T »x}{»x-v} .  
	\]
	This proves the pseudomonotonicity of $ T $.
\endproof

In order to study systems, we consider the following continuity property.

\begin{definition}[{\rm Sequential solutional continuity}]\label{..}
	Suppose that $ X $ and $»Z$ are two topological spaces,
		$ Y $ is an arbitrary set and that $ T : X \times Y  \rightarrow  »Z$.
	We say that $ T $ is sequentially solutionally continuous in $»x\in X $ and $»z\in»Z$ if the equation
	$ T (»x,»y)=»z$ has a unique solution $»y\in Y $, and if for every sequence $( x_n )_{n\in{{\mathbb N}}}$ converging
	to $»x$ in $ X $ holds
	\[
		T ( x_n ,»y) \>\mathop{\rightarrow }\>  »z \enspace \enspace \mbox{in }»Z.
	\]  
	Furthermore, $ T $ is said to be sequentially solutionally continuous in $»z\in»Z$ if $ T $ is so in
	$»x$ and $»z$ for every $»x\in X $.
\end{definition}

Next, assumptions are given that guarantee the pseudomonotonicity of the operator $ S $ from the introduction.
We suppose the uniformly strong monotonicity and the sequential solutional continuity of $»\xwt{ B }$ as well as Lipschitz conditions.
The assumptions {\rm(A3.2)} and {\rm(A3.3)} can be seen as a counterpart to conditions (S2)--(S4)
	used in the definition of semimonotone operators.


\begin{definition}[{\rm Assumptions {\rm(A1)}, {\rm(A2)} and {\rm(A3)}}]\label{a2:Semi}{\parskip3pt%
	\def\xa#1#2{\rucky{2cm}{\hfil{\rm #1}\hfil}{#2}}%
	\def\xb#1#2#3#4{\pboxc{2cm}{\rm #1}\hskip10mm\pboxl{7.5cm}{$#2$}\pboxc{1cm}{$#3$}\pboxl{2cm}{$#4$} }%
	\def\xc#1#2#3{\pboxc{2cm}{\rm #1}\hbox to 13.5cm{#2\hss$#3$}}%
	\def\xd#1#2#3#4{\pboxl{2.5cm}{#1}\pboxc{6cm}{$#2$\pboxc{7mm}{$#3$}$#4$}}%

	We say the {\rm(A1)} is fulfilled if the following conditions are met:
		\xa{{\rm(A1.1)}}{$X,Y$ are real, reflexive Banach spaces and
				and $ y_0^*\in Y ^*$,}
		\xa{{\rm(A1.2)}}{$ A : X \times Y  \rightarrow   X ^*$ and $ B : X \times Y  \rightarrow   Y ^*$ together with
				$\xwt{ A }: X \times X \times Y   \rightarrow   X ^*$ and $\xwt{ B }: X \times X \times Y   \rightarrow   Y ^*$ are mappings such that
				$ A (»x,»y)=\xwt{ A }(»x,»x,»y), \enspace  B (»x,»y)=\xwt{ B }(»x,»x,»y)$
				for all $(»x,»y)\in X \times Y  $.
				}
		\xa{{\rm(A1.3)}}{The mapping $»y\mapsto \xwt{ B }( x_1 , x_2 ,»y)$ from $ Y $ into $ Y ^*$ is strongly monotone uniformly in
				$( x_1 , x_2 )\in X \times X $, i.e.\ there exists an $\alpha_B >0$ such that
				\[   \iKlammerB{\xwt{ B }( x_1 , x_2 , y_1 )-\xwt{ B }( x_1 , x_2 , y_2 )}{ y_1 - y_2 }  \iSub{ Y }  \>\ge\>    \alpha_B \, ||\iBlock{ y_1 - y_2 }||  _{ Y } ^2  \]
				for all $ x_1 , x_2 \in X  $ and $ y_1 , y_2 \in Y $.
				Furthermore, $»y\mapsto \xwt{ B }( x_1 , x_2 ,»y)$ is radially continuous for every tuple $( x_1 , x_2 )\in X  \times X  $.}
	If furthermore there are constants $\beta_A ,\beta_B \ge0$ and $\alpha_A >0$ such that\\[12pt]
		\xb{{\rm(A2.1)}}{	 \iKlammerB{\xwt{ A }( x_1 , x_2 ,»y)		- \xwt{ A }( x_2 , x_2 ,»y)}{ x_1 - x_2 }  \iSub{ X }	}{\ge}{  \alpha_A \, ||\iBlock{ x_1 - x_2 }||  _{ X } ^2, }   \\[0.5ex]
		\xb{{\rm(A2.2)}}{	 ||\iBlock{\xwt{ A }( x_1 , x_2 , y_1 )	- \xwt{ A }( x_1 , x_2 , y_2 )}||  _{ X ^*} 			}{\le}{  \beta_A \, ||\iBlock{ y_1 - y_2 }||  _{ Y } ,   }   \\[0.5ex]
		\xb{{\rm(A2.3)}}{	 ||\iBlock{\xwt{ B }( x_1 , x_2 ,»y)	- \xwt{ B }( x_2 , x_2 ,»y)}||  _{ Y ^*} 			}{\le}{  \beta_B \, ||\iBlock{ x_1 - x_2 }||  _{ X }     }   \\[12pt]
	hold for all $ x_1 , x_2 \in X  $ and $»y\in Y $, then we say that {\rm(A2)} is satisfied.
	Finally, {\rm(A3)} is fulfill if {\rm(A2)} \!and the following conditions are satisfied \\[2ex]
		\xc{{\rm(A3.1)}}{the mapping $(»x,»y)\mapsto \xwt{ B }(»x_0 ,»x,»y)$ from $ X _\omega \times Y $ into $ Y ^*$                          }{              }  \\[0.5mm]
		\xc{    }{is sequentially solutionally continuous in $»x=»x_0 $ and $ y_0^*$                                                            }{ \forall»x_0 \in X  , }  \\[3mm]
		\xc{{\rm(A3.2)}}{if $ x_n \mathop{\rightharpoonup}»x,\enspace  y_n \mathop{\rightarrow }»y$ and $\mathop{\overline {\lim}}\limits_{n\mathop{\rightarrow }\infty}  \iKlammerB{\xwt{ A }(»x, x_n , y_n )}{ x_n -»x}  \iSub{ X }\le0$  }{              }  \\[0.5ex]
		\xc{    }{\xd{then it holds }{ \xwt{ A }(»x_0 , x_n , y_n )              }{ \mathop{\rightharpoonup} }{ \xwt{ A }(»x_0 ,»x,»y)\enspace  }      }{              }  \\[0.5ex]
		\xc{    }{\xd{and           }{  \iKlammerB{\xwt{ A }(»x_0 , x_n , y_n )}{ x_n -»x}  \iSub{ X } }{ \mathop{\rightarrow }  }{ 0                 }      }{ \forall»x_0 \in X  , }  \\[3mm]
		\xc{{\rm(A3.3)}}{the mapping $»x\mapsto \xwt{ A }(»x,»x_0 ,»y)$ is radially continuous in $»x_0 $                                         }{ \forall»x_0 \in X  ,\>»y\in Y , } \\[3mm]%
		\xc{{\rm(A3.4)}}{ $\alpha_A \,\alpha_B  \>\ge\>   \beta_A \,\beta_B $                                                  }{              }\\[2ex]
  for all sequences $( x_n )_{n\in{{\mathbb N}}}$ in $ X $ and $( y_n )_{n\in{{\mathbb N}}}$ in $ Y $.
}\end{definition}

Particularly, for every $»y\in Y $ the mapping $»x\mapsto  A (»x,»y)$ is semimonotone.
Moreover, since $\xwt{ B }_{( x_1 , x_2 )}: Y \rightarrow   Y ^*$ is strongly monotone and radially continuous,
	the equation $\xwt{ B }_{( x_1 , x_2 )}»y= y_0^*$ has a unique solution $»y\in Y $ for every $ x_1 , x_2 \in X  $.
The corresponding solution operator and its composition with $\xwt{ A }$ are denoted by
	$»\xwt{ R }$ and $»\xwt{ S }$, respectively.

\begin{definition}[{\rm Operators $\xwt{ R }$ and $\xwt{ S }$}]\label{d2:RS}{\def\xa#1{\pboxl{3cm}{$#1$}}%
	\def\xb#1{\pboxl{8cm}{$#1$}}%
	\def\xc#1#2{\xa{#1} && \xb{#2}}%
	\def\xd#1#2{\xa{#1} \hskip7mm \xb{#2}}%
	Assume {\rm(A1)} and $ x_1 , x_2 \in X  $.
	The bijectivity of $\xwt{ B }_{( x_1 , x_2 )}$ allows us to define the operators
		$\xwt{ R }$ and $\xwt{ S }$ on $ X \times X  $ into $ Y $ respectively $ X ^*$ as
	\begin{eqnarray*} 
			\xc{ \xwt{ R }: X \times X  \rightarrow   Y ,  }{ \xwt{ R }( x_1 , x_2 ):=\xwt{ B }_{( x_1 , y_2 )}^{\>-1}\,( y_0^*). } \\[0.5ex]
			\xc{ \xwt{ S }: X \times X  \rightarrow   X ^*, }{ \xwt{ S }( x_1 , x_2 ):=\xwt{ A }( x_1 , x_2 ,\xwt{ R }( x_1 , x_2 )).     }  
	\end{eqnarray*}%
}\end{definition}

The following lemma provides simple Lipschitz and monotonicity properties of $\xwt{ R }$ and $\xwt{ S }$.

\begin{lemma}\label{p1:LipMon}
	If the {\rm(A1)} is fulfilled, then
	\[
		 ||\iBlock{ \xwt{ R } z_1 -\xwt{ R } z_2  }||  _{ Y }   \enspace \le\enspace    \f1{\alpha_B } \> ||\iBlock{ \xwt{ B }(»z,\xwt{ R } z_1 )-\xwt{ B }(»z,\xwt{ R } z_2 ) }||  _{ Y ^*} 
	\]
	holds for pairs $»z, z_1 , z_2 \in X  \times X  $.
	If {\rm(A2)} is satisfied, then for all $ x_1 , x_2 \in X  $
	\begin{eqnarray*} 
		 ||\iBlock{\xwt{ R }( x_1 , x_2 )-\xwt{ R }( x_2 , x_2 )}||  _{ Y }   & \>\le\>   & \f{\beta_B }{\alpha_B }\,  ||\iBlock{ x_1 - x_2 }||  _{ X } ,     \\[0.5ex]
		 \iKlammerB{\xwt{ S }( x_1 , x_2 )-\xwt{ S }( x_2 , x_2 )}{ x_1 - x_2 }  \iSub{ X }  & \>\ge\>   & \f{\alpha_A \alpha_B -\beta_A \beta_B }{\alpha_B }\, ||\iBlock{ x_1 - x_2 }||  _{ X } ^2.
	\end{eqnarray*}%
\end{lemma}

{\em Proof}. {
	For $»z, z_1 , z_2 \in X  \times X  $ {\rm(A1.3)} implies
		\begin{eqnarray*} 
			 ||\iBlock{ \xwt{ R } z_1 -\xwt{ R } z_2  }||  _{ Y } ^2   & \le &  \f{1}{\alpha_B }\>   \iKlammerB{ \xwt{ B }(»z,\xwt{ R } z_1 )-\xwt{ B }(»z,\xwt{ R } z_2 )}{\xwt{ R } z_1 -\xwt{ R } z_2 }  \iSub{ Y } \\[0.5ex]
					& \le &   \f{1}{\alpha_B }\>   ||\iBlock{ \xwt{ B }(»z,\xwt{ R } z_1 )-\xwt{ B }(»z,\xwt{ R } z_2 )}||  _{ Y ^*}  \,  ||\iBlock{\xwt{ R } z_1 -\xwt{ R } z_2 }||  _{ Y } 
		\end{eqnarray*}%
		and hence the first inequality.
	Assuming {\rm(A2)} and $ x_1 , x_2 \in X $ and with $ »z_i :=( »x_i , x_2 )$, from the definition of $\xwt{ R }$ we have
		$\xwt{ B }( »z_i ,\xwt{ R } »z_i )= y_0^*$ and therefore by the first inequality that
		\[
			 ||\iBlock{ \xwt{ R } z_1 -\xwt{ R } z_2  }||  _{ Y }  \enspace \le\enspace   \f1{\alpha_B }\>  ||\iBlock{\xwt{ B }( z_2 ,\xwt{ R } z_2 ) - \xwt{ B }( z_1 ,\xwt{ R } z_2 )}||  _{ Y ^*}   \enspace \le\enspace   \f{\beta_B }{\alpha_B }\>  ||\iBlock{ x_1 - x_2 }||  _{ X } ,
		\]
		which is the second inequality.
	Together with {\rm(A2.1)} and {\rm(A2.2)}, this yields the estimation
	\begin{eqnarray*} 
		&& \hb     \iKlammerB{ \xwt{ S }( x_1 , x_2 )-\xwt{ S }( x_2 , x_2 ) }{ x_1 - x_2 }  \iSub{ X }   \\[0.5ex]
		&  =  &    \iKlammerB{\xwt{ A }( x_1 , x_2 ,\xwt{ R }( x_1 , x_2 ))-\xwt{ A }( x_2 , x_2 ,\xwt{ R }( x_1 , x_2 )) }{ x_1 - x_2 }  \iSub{ Y }   \\[0.5ex]
		&     &   +\enspace   \iKlammerB{\xwt{ A }( x_2 , x_2 ,\xwt{ R }( x_1 , x_2 ))-\xwt{ A }( x_2 , x_2 ,\xwt{ R }( x_2 , x_2 )) }{ x_1 - x_2 }  \iSub{ Y }   \\[0.5ex]
		& \ge &   \alpha_A \,  ||\iBlock{ x_1 - x_2 }||  _{ X } ^2   »-  \beta_A \,  ||\iBlock{\xwt{ R }( x_1 , x_2 )-\xwt{ R }( x_2 , x_2 )}||  _{ Y }   ||\iBlock{ x_1 - x_2 }||  _{ X }   \\[0.5ex]
		& \ge &   \f{\alpha_A \alpha_B -\beta_A \beta_B }{\alpha_B } \>\>  ||\iBlock{ x_1 - x_2 }||  _{ X } ^2,
	\end{eqnarray*}%
	and finishes the proof.
}\endproof


The following lemma is crucial in order to prove the pseudomonotonicity of $\xwt{ B }$.

\begin{lemma}\label{p1:RRcont}{\def\xa#1#2{\pboxr{1cm}{$#1$} \enspace \mapsto \enspace  \pboxl{2cm}{$#2$}}%
	Suppose {\rm(A1)}, $ x_0 \in X  $ and that $ X _T $ denotes $ X $ equipped with some topology $T$.
	If the mapping $ \xwt{ B }_{ x_0 } : X _T  \times Y \rightarrow   Y ^*$ is sequentially solutionally continuous in $ x_0 $ and $ y_0^*$, 
		then $ \xwt{ R }_{ x_0 } : X _T  \rightarrow   Y $ is continuous in $ x_0 $.
}\end{lemma}

{\em Proof}. 
	Let $( x_n )_{n\in{{\mathbb N}}}$ be a sequence in $ X  $ converging to $ x_0 $ with respect to $ X _T  $. Lemma~\ref{p1:LipMon}
	provides the estimation
	\[
			 ||\iBlock{\xwt{ R }( x_0 , x_n )-\xwt{ R }( x_0 , x_0 )}||  _{ Y } 
		\enspace \le\enspace  	\f1{\alpha_B }\,  ||\iBlock{\xwt{ B }( x_0 , x_n ,\xwt{ R }( x_0 , x_n ))-\xwt{ B }( x_0 , x_n ,\xwt{ R }( x_0 , x_0 ))}||  _{ Y ^*} 
	\]
	for all $»n\in{{\mathbb N}}$.
	Furthermore, $\xwt{ B }( x_0 , x_n ,\xwt{ R }( x_0 , x_n ))= y_0^*=\xwt{ B }( x_0 , x_0 ,\xwt{ R }( x_0 , x_0 ))$ from the definition of $\xwt{ R }$.
	Hence, the sequential solutional continuity of $\xwt{ B }$ implies that $\xwt{ R }( x_0 , x_n )$ converges strongly to $\xwt{ R }( x_0 , x_0 )$ in $ Y $.
\endproof

The next theorem provides the semimonotonicity and hence the pseudomonotonicity of $ S $.

\begin{theorem}[{\rm Semimonotone Reduction}]\label{t1:SemReduct}{\def\xa#1#2{\pboxl{6mm}{{\rm #1}}{#2}}%
	Suppose {\rm(A3)} and $ y_0^*\in X ^*$, and let $»F$ be the mapping $( A , B )»: X  \times Y \rightarrow   X ^*\times Y ^*$.
	Then the operators $ R $ and $ S $ of Definition~\ref{d2:RS} satisfy the following statements:\\[0.5ex]
	\xa{i)}{$	»F(»x,»y) »= ( x_0^*, y_0^*) \enspace \enspace \Lral\enspace \enspace  »y= R »x$ \enspace  and \enspace $ S »x »= x_0^*$,			} \\[0.5ex]
	\xa{ii)}{$	S : X  \rightarrow   X ^* $ \enspace  is semimonotone with the semimonotone representative $\xwt{ S }: X \times X \rightarrow   X ^*$.		}\\[0.5ex]
}\end{theorem}

{\em Proof}. 
	Part i) follows from {\rm(A1)}.
	To prove ii) we show that the operator $\xwt{ S }$ satisfies the conditions (S1)--(S4). 

	The condition (S1) immediately follows from Lemma~\ref{p1:LipMon} combined with condition~{\rm(A3.4)}.
	To show (S2) and (S3), let us consider an $ S $--pas $( x_n )_{n\in{{\mathbb N}}}$ which weakly converges
	to $»x\in X $. The monotonicity property~(S1) of $\xwt{ S }$ yields
	\begin{equation} \label{«cs1}
		 \iKlammerB{\xwt{ S }(»x, x_n )}{ x_n -»x}  \iSub{ X }  \>\le\>    \iKlammerB{\xwt{ S }( x_n , x_n )}{ x_n -»x}  \iSub{ X }  »=  \iKlammerB{ S x_n }{ x_n -»x}  \iSub{ X }.
	\end{equation} %
	Passing to the limit superior on both sides and using the $ S $--pas property of $( x_n )$ shows
	\[
		\mathop{\overline {\lim}}_{n\mathop{\rightarrow }\infty}  \iKlammerB{\xwt{ A }(»x, x_n ,\xwt{ R }(»x, x_n )}{ x_n -»x}  \iSub{ X }  »= \mathop{\overline {\lim}}_{n\mathop{\rightarrow }\infty}  \iKlammerB{\xwt{ S }(»x, x_n )}{ x_n -»x}  \iSub{ X }  \>\le\>   0.
	\]
	The sequential solutional continuity~{\rm(A3.1)} together with Proposition~\ref{p1:RRcont} yields
	\begin{equation} \label{«cs2}
		\xwt{ R }(»x, x_n ) \mathop{\relbar\joinrell\rightarrow} \xwt{ R }(»x,»x)   \enspace \enspace \mbox{in\enspace } Y.
	\end{equation} %
	Thus, we can apply {\rm(A3.2)} in order to obtain
	\begin{eqnarray} 
		& \xwt{ S }(»x, x_n ) »= \xwt{ A }(»x, x_n ,\xwt{ R }(»x, x_n )) \mathop{\relbar\joinrell\rightharpoonup} \xwt{ A }(»x,»x,\xwt{ R }(»x,»x)) »= \xwt{ S }(»x,»x), &    \nonumber\\[4pt]
		& \lim\limits_{n\mathop{\rightarrow }\infty}  \iKlammerB{\xwt{ S }(»x, x_n )}{ x_n -»x}  \iSub{ X }  »=  \lim\limits_{n\mathop{\rightarrow }\infty}  \iKlammerB{\xwt{ A }(»x, x_n , R (»x, x_n )}{ x_n -»x}  \iSub{ X }  »=  0.&  \label{«cs3}
	\end{eqnarray}%
	These are the properties (S2) and (S3).
	Finally, it is easy to check that the radial continuity {\rm(A3.3)} in combination with the
	Lipschitz properties {\rm(A2.2)} and \nix(\ref{«cs2}) imply that the mapping
	\[
		»x \>\mapsto \> \xwt{ S }(»x,»x_0 ) »= \xwt{ A }(»x,»x_0 ,\xwt{ R }(»x,»x_0 ))
	\]
	is radially continuous in $»x_0 \in X  $. This shows (S4) and therefore completes the proof.
\endproof

\begin{remark}\label{..}
	Assume $\alpha_A \alpha_B >\beta_A \beta_B $. Then, by Lemma~\ref{p1:LipMon} the operator $\xwt{ S }$ satisfies a
		strong monotonicity condition in the first argument.
	Hence, \nix(\ref{«cs1}) can be strengthened to
	\[
			 \iKlammerB{\xwt{ S }(»x, x_n )}{ x_n -»x}  \iSub{ X } + c\, ||\iBlock{ x_n -»x}||  _{ X } ^2  \>\le\>    \iKlammerB{ S x_n }{ x_n -»x}  \iSub{ X }
	\]
	with $c:=\f{1}{\alpha_B }(\alpha_A \alpha_B -\beta_A \beta_B )>0$.
	This together with the $ S $--pas condition on $( x_n )$ and the convergence
		$ \iKlammerB{\xwt{ A }(»x_0 , x_n , y_n )}{ x_n -»x}  \iSub{ X } \mathop{\rightarrow }  0  $ \>from~{\rm(A3.2)}
		shows that $( x_n )$ even converges strongly to $»x$.
	Consequently, if $\alpha_A \alpha_B  > \beta_A \beta_B $, we can relax {\rm(A3)} by requiring the desired convergence properties
		in {\rm(A3.2)}	only if $ x_n \mathop{\rightarrow }»x$, $ y_n \mathop{\rightarrow }»y$ and $\mathop{\overline {\lim}}\limits_{n\mathop{\rightarrow }\infty}  \iKlammerB{\xwt{ A }(»x, x_n , y_n )}{ x_n -»x}  \iSub{ X }\le0.$
	
%
\end{remark}

The final proposition of this section ensures the demicontinuity of $ S $.

\begin{proposition}\label{p1:Demi}
	Assume {\rm(A2)} and suppose for every $ x_0 \in X  $ and $»y\in Y $ that
	{\def\xa#1{»x\>\mapsto \>\pboxl{1.3cm}{$#1$}}%
	\begin{eqnarray*} 
			\xa{\xwt{ R }( x_0 ,»x)}  && \mbox{is continuous,} \\[0.5ex]
			\xa{ A (»x,»y)}    && \mbox{is demicontinuous.}
	\end{eqnarray*}}%
	Then $ S : X  \rightarrow   X ^*$ is demicontinuous.
\end{proposition}

{\em Proof}. 
	Assume that $ x_n \mathop{\rightarrow }»x$. Lemma~\ref{p1:LipMon} and the continuity of $ R $ imply
	\begin{eqnarray*} 
		&& \hb
		\lim_{n\mathop{\rightarrow }\infty}  ||\iBlock{ R x_n - R »x}||  _{ X }   \\[0.5ex]
		& \le &   \lim_{n\mathop{\rightarrow }\infty}  ||\iBlock{\xwt{ R }( x_n , x_n )-\xwt{ R }(»x, x_n )}||  _{ X }   »+ \lim_{n\mathop{\rightarrow }\infty}  ||\iBlock{\xwt{ R }(»x, x_n )-\xwt{ R }(»x,»x)}||  _{ X }  \\[0.5ex]
		&  =  &   0.
	\end{eqnarray*}%
	Thus, condition~{\rm(A2.2)} in combination with the demicontinuity of $ A $ yields
	\begin{eqnarray*} 
		&& \hb  \lim_{n\mathop{\rightarrow }\infty}  \iKlammerB{ S x_n - S »x}{»v}  \iSub{ X }  					\\[0.5ex]
		& = &   \lim_{n\mathop{\rightarrow }\infty}  \iKlammerB{ A ( x_n , R x_n )- A ( x_n , R »x)}{»v}  \iSub{ X }		
				+\lim_{n\mathop{\rightarrow }\infty}  \iKlammerB{ A ( x_n , R »x)- A (»x, R »x)}{»v}  \iSub{ X } 		\\[0.5ex]
		& = &   0
	\end{eqnarray*}%
	for every $»v\in X $. This finishes the proof.
\endproof

\begin{remark}\label{«csr1}
	The continuity assumption on $\xwt{ R }$ is fulfilled, for instance, if {\rm(A3.1)} holds (cf. Lemma~\ref{p1:RRcont}).
	Moreover, if $ A $ is continuous in the first argument, then $ S $ is continuous.
\end{remark}



\section{Abstract evolution equations}
\label{s3:}

Before turning to a special application of Theorem~\ref{t1:SemReduct} in the next section,
	we present some elements of the framework of Gr\"oger~\cite{konni} for evolution equations
	allowing to include compositions with certain linear operators under the time derivative.
Well-known embedding theorems and results on existence, uniqueness or the continuous dependence on
	the data hold also within this framework.
For further details and proofs we refer to~\cite{konni,jens,doni}.

Throughout this section we suppose the following.

\begin{assumption}\label{a2:Spaces}
	Let $ V $ be a reflexive Banach space such that $ V $ and $ V ^*$ are strictly convex,
		$ H $ a Hilbert space and $ K \in L( V ; H )$ be an operator 
		having dense image $ K ( V )$ in $ H $.
		The operator $ E \in L( V ; V ^*)$ is given by $ E := K ^* J_H K $ 
		«($ J_H $ is the duality mapping of~$ H $«).
	Moreover, suppose that ${\cal T}=\mathop{]\kern1pt  0 ,T\kern1pt[}$ with $T>0$ and $1<p,p'<\infty$ with
		$\f1p+\f1{p'}=1$ and $p\ge2$.
\end{assumption}

\begin{remark}\label{r3:E}\nix\\
	\ifnum\the\remarkcc=0\else \\\fi\advance\remarkcc by 1{\bf \the\remarkcc.}\enspace  %
				The operator $ E \in L( V ; V ^*)$ is positive and symmetric
				(i.e.	$	 \iKlammerB{ E »u}{»u}  \ge 0,$
						$	 \iKlammerB{ E »u}{»v} = \iKlammerB{ E »v}{»u} 	$
						$	\forall »u,»v\in V )$.
				Conversely, given any positive and symmetric operator $ E \in L( V ; V ^*)$
					we can choose $H$ as the completion of the
					pre-Hilbert space $ V /\ker E $ with inner product $ \iKlammerA{»u}{»v} := \iKlammerB{ E »u}{»v} $
					and $ K »u:=[»u]$ in order to satisfy~Assumption~\ref{a2:Spaces}.\ifnum\the\remarkcc=0\else \\\fi\advance\remarkcc by 1{\bf \the\remarkcc.}\enspace  %
				If $ K $ is injective, it is a bijection from $ V $ onto $ K ( V )\subset H $.
				Therefore, $ V $ and $ H $ can be regarded as a usual evolution triple by
					identifying $ V $ with $ K ( V )$,
					equipping $ K ( V )$ with the norm
					$ ||\iBlock{x}||  _{ K ( V )} := ||\iBlock{ K ^{-1} x}||  _{ V } $ and 
					considering the $ K ( V )\hookrightarrow  H \cong H ^*\hookrightarrow ( K ( V ))^*$.
				We use this identification of\/ $ V $ with $ K ( V )$ even if
					$ V $ is a subset of $ H $ itself, cf.~Section~\ref{s4:}.
				\end{remark}

Corresponding to these spaces and operators 
	we define ${\cal  V }:=L^2({\cal T}; V )$ and ${\cal  H }:=L^2({\cal T}; H )$ with standard norms
	and identity ${\cal  V }^*$ with $L^2({\cal T}; V ^*)$
	(which we can do since $ V $ is reflexive and therefore possesses the Radon-Nikod\'ym property, cf. \cite{uhl}).
Moreover, we set
	$
		({\cal E}»u)(t):= E »u(t)$ and $
		({\cal  K }»u)(t):= K »u(t)
	$
	in order to obtain ${\cal E}\in L({\cal  V };{\cal  V }^*)$ and ${\cal  K }\in L({\cal  V };{\cal  H })$.
The space ${\cal W}$ is the space of all $»u\in{\cal  V }$ such that ${\cal E}»u\in{\cal  V }^*$ possesses a
	weak time derivative which again belongs to ${\cal  V }^*$:
	\begin{eqnarray*} 
			{\cal W}			:=	\{»u\in{\cal  V } \>|\>  {\cal E}»u \mbox{ has a weak derivative } ({\cal E}»u)'\in{\cal  V }^*\},	\enspace 
			 ||\iBlock{»u}||  _{{\cal W}} 	:=	( ||\iBlock{»u}||  _{{\cal  V }} ^2 »+  ||\iBlock{({\cal E}»u)'}||  _{{\cal  V }^*} ^2)^{1/2}.								
	\end{eqnarray*}
Furthermore, we define the linear operator ${\cal  L }\subset{\cal  V }\times{\cal  V }^*$ by
	\[
		D({\cal  L })	:=	{\cal W}\subset{\cal  V },			\enspace \enspace 
		{\cal  L }»u:=({\cal E}»u)'\in{\cal  V }^*
	\]
	and ${\cal I}\in L({\cal W};{\cal  V })$ as the identity ${\cal I}:=\mathop{\rm Id}\rulenix_{{\cal W}\mathop{\rightarrow }{\cal  V }}$
	regarded as a mapping from ${\cal W}$ into ${\cal  V }$.
For these spaces we obtain the following density result and a formula of integration by parts.

\begin{proposition}\label{lab.prop.fevol.banach}
	The space ${\cal W}$ is a reflexive Banach space and
		$\{»u|_{{\cal T}} »: »u\in C^\infty_c({{\mathbb R}}; V )\}$ is a dense subspace.
\end{proposition}

\begin{proposition}\label{p3:PartInt}
	The operator ${\cal  K }$ maps ${\cal W}$ continuously into the space $C(\overline {{\cal T}}; H )$, meaning that every class of
		equivalent functions in ${\cal  K }({\cal W})\subset L^p({\cal T}; H )$ possesses a representative that is
		continuous from ${\cal T}$ into $ H $ with continuous extension onto $\overline {{\cal T}}$.
	Furthermore, in this sense the formulas hold for all $»u,»v\in{\cal W}$ and $ t_1 , t_2 \in\overline {{\cal T}}$
	\begin{eqnarray*} 
		&				 \iKlammerA{({\cal  K }»u)( t_2 )}{({\cal  K }»v)( t_2 )}  \iSub{ H } »-  \iKlammerA{({\cal  K }»u)( t_1 )}{({\cal  K }»v)( t_1 )}  \iSub{ H }			\hskip5cm				&\\
		&\hskip5cm			\>=\> 	\int_{ t_1 }^{ t_2 } \big[  \iKlammerB{({\cal E}»u)'(t)}{»v(t)}  \iSub{ V } »+  \iKlammerB{({\cal E}»v)'(t)}{»u(t)}  \iSub{ V }  \big]  \,dt,		&\\[4pt]
		&				 ||\iBlock{({\cal  K }»u)( t_2 )}||  _{ H } ^2 »-  ||\iBlock{({\cal  K }»u)( t_1 )}||  _{ H } ^2 
								\>=\>  2\int_{ t_1 }^{ t_2 }  \iKlammerB{({\cal E}»u)'(»t)}{»u(»t)}  \iSub{ V }\,dt.					&
	\end{eqnarray*}%
\end{proposition}

In order to incorporate the treatment of initial data of evolution equations directly into
	the operators and the spaces let us consider
	\[
			\widehat {{\cal W}} «:= {\cal W}\times H ,	\enspace \enspace \enspace 
			\widehat {{\cal  V }} «:= {\cal  V }\times H ,
	\]
	with the product norm $ ||\iBlock{(x,y)}||  _{X\times Y} :=( ||\iBlock{x}||  _{X} ^2+ ||\iBlock{y}||  _{Y} ^2)^{1/2}$ on $X\times Y$ for
		two normed vector spaces $X$ and $Y$.
{\def\xa#1#2#3{&\pboxl{1.3cm}{$#1$}\pboxl{6cm}{$#2$}\pboxl{4cm}{$#3$}&}%
To a given $»h\in H $ the (single-valued) operators $\widehat {{\cal  L }}\subset \widehat {{\cal  V }}\times \widehat {{\cal  V }}^*$ and
	${\cal  L }_{»»h}\in {\cal  V }\times {\cal  V }^* $ are defined by
	\begin{eqnarray*} 
			\xa{ D(\widehat {{\cal  L }})  }{ :=\{(»u,({\cal  K }»u)(0)) \>|\>  »u\in{\cal W}\},	}{  \rulenix \widehat {{\cal  L }}(»u,»h):=({\cal  L }»u, J_H »h),  } \\[0.5ex]
			\xa{ D({\cal  L }_{»»h}) }{ :=\{»u\in{\cal W} »: ({\cal  K }»u)(0)=»h\},	}{  {\cal  L }_{»»h}:={\cal  L }|_{D({\cal  L }_{»»h})}.  }
	\end{eqnarray*}%
}%

A fundamental result is the maximal monotonicity of $\widehat {{\cal  L }}$.

\begin{proposition}\label{}
	The operator $\widehat {{\cal  L }}\in \widehat {{\cal  V }}\times \widehat {{\cal  V }}^* $ is a linear, maximal monotone operator.
\end{proposition}

\begin{corollary}\label{}
	For every $»h\in H $, the operator ${\cal  L }_{»»h}\in {\cal  V }\times {\cal  V }^* $ is maximal monotone.
\end{corollary}

{\em Proof}. 
	By \cite[Theorem~32.F]{zeid2b} a monotone mapping $ T \subset X \times X ^*$ on a reflexive Banach space $ X $
		with $ X $ and $ X ^*$ being strictly convex
		is maximal monotone if and only if $ T + J_{ X }$ is surjective.
	Therefore, let an arbitrary $»u^*\in{\cal  V }^*$ be given. 
	Applied to $\widehat {{\cal  L }}$, the theorem in question  shows the existence of a
		$\widehat {»u}=(»u, h_1 )\in\widehat {{\cal  V }}$ such that $(\widehat {{\cal  L }}+»J_{\widehat {{\cal  V }}})\widehat {»u}=(»u^*,2 J_H »h)$.	
	Since $\widehat {»u}\in D(\widehat {{\cal  L }})$, it follows that $ h_1 =({\cal  K }»u)(0)$.
	Moreover, it is easy to check that $ J_{\widehat {{\cal  V }}}(»v,»g)=( J_{{\cal  V }}»v, J_H »g)$.
	This implies $({\cal  L }+»J_{{\cal  V }})»u=»u^*$ and $2 J_H ({\cal  K }»u)(0)=2 J_H »h$.
	Consequently, we conclude that $»u\in D({\cal  L }_{»»h})$ and $({\cal  L }_{»»h}+»J_{{\cal  V }})»u=»u^*$.
\endproof

The following theorems provides conditions that ensure the solvability of evolution inclusions.

\begin{theorem}\label{t2:Exist}
	Suppose Assumption~\ref{a2:Spaces} and $(»f, h )\in{\cal  V }^*\times H $.
	Let ${\cal A}:{\cal  V }\rightarrow  {\cal  V }^*$ be bounded, maximal monotone and
		${\cal B}:{\cal  V }\rightarrow  {\cal  V }^*$ be bounded, demicontinuous, coercive with respect to a $ w_0 \in{\cal W}$ such that
			${\cal B}$ is pseu\discretionary{-}{}{}domono\discretionary{-}{}{}tone with respect to ${\cal  L }_{0}$.
	Then there exists a $»u\in{\cal W}$ with
	\[
			({\cal  L }+{\cal A}+{\cal B})»u »= »f,\enspace \enspace \enspace  ({\cal  K }»u)(0)= h .   
	\]
\end{theorem}

{\em Proof}. 
	1. Let us choose a $»w\in D({\cal  L }_{- h })$ (note that ${\cal  L }_{- h }\neq\emp$ since it is maximal monotone).
	First, we show that there exists a $»v\in D({\cal  L }_{0})$ with
	\[
			({\cal  L }_{0}+  \mathfrak{T}_{»w}{\cal A}+  \mathfrak{T}_{»w}{\cal B})»v »= »f+{\cal  L }»w.
	\]
	Since $\mathfrak{T}_{»w}{\cal A}$ is bounded and maximal monotone, it is pseu\discretionary{-}{}{}domono\discretionary{-}{}{}tone and demicontinuous (cf.~\cite[Lemma~1.3, p.~66]{ggz}).
	Particularly, $\mathfrak{T}_{»w}{\cal A}$ is pseu\discretionary{-}{}{}domono\discretionary{-}{}{}tone with respect to ${\cal  L }_{0}$ as well as $\mathfrak{T}_{»w}{\cal B}$ (cf.~Proposition~\ref{p2:PsmSum}).
	Consequently, $\mathfrak{T}_{»w}{\cal A}+\mathfrak{T}_{»w}{\cal B}$ is pseu\discretionary{-}{}{}domono\discretionary{-}{}{}tone with respect to ${\cal  L }_{0}$, demicontinuous and
		coercive with respect to $ w_0 +»w$.
	An existence result by Lions~\cite[Theorem~1.1, p.~316]{lions} guarantees that
		there is a $»v\in D({\cal  L }_{0})$ with
		$({\cal  L }_{0}+\mathfrak{T}_{»w}{\cal A}+\mathfrak{T}_{»w}{\cal B})»v\ni»f+{\cal  L }»w$.

	2. Setting $»u:=»v-»w$ it holds $»u\in D({\cal  L }_{»»h})$ and ${\cal  L }_{»»h}»u = {\cal  L }_{0}»v - {\cal  L }»w$.
	This implies
	\[
			({\cal  L }_{»»h}+{\cal A}+{\cal B})»u «= ({\cal  L }_{0}+\mathfrak{T}_{»w}{\cal A}+\mathfrak{T}_{»w}{\cal B})»v - {\cal  L }»w  «=  »f
	\]
	and therefore completes the proof.
\endproof

\begin{remark}\label{}
	The theorem by Lions applied in our proof assumes coercivity of the pseu\discretionary{-}{}{}domono\discretionary{-}{}{}tone operator with
		respect to $0$, but it can be generalized to the case of
		coercivity with respect to an arbitrary element in ${\cal W}$ without any difficulties
		(cf.~also \cite[Theorem~2.6.1]{doni}).
\end{remark}


\section{Application to a model of phase separation}
\label{s4:}

In this section we show how the results of the previous »sections can be applied to prove
	the existence of solutions to coupled elliptic-parabolic systems.
In order to demonstrate the ability of these techniques and the generality of Gr\"oger's framework
	we consider a parabolic equation of fourth order in space of Cahn--Hilliard type
	which is coupled to a elliptic equation modeling a quasi-steady mechanical equilibrium
	for each point in time.
The given system is highly nonlinear and both parts are strongly coupled.
This generality imposes a restriction -- solving the elliptic part and inserting the solution into the
	parabolic part we use Theorem~\ref{t1:SemReduct} to ensure the pseudomonotonicity of the resulting operator.
Therefore, we require Assumption~\ref{a2:Semi} to hold which means that we have to restrict
	the influence both parts of the system may exert on each other.
This is necessary since changes in lower order terms of one part may effect
	higher order terms in the other and
	the reduced equation has to be monotone in the leading order terms.
Nevertheless, no other existence results for this very general system
	seem to be known yet.

Together with initial and boundary conditions, our systems reads as follows
{\def\xa#1#2#3{&\pboxc{0cm}{\pboxc{0.87\hsize}{$#1$}\pboxl{0.04\hsize}{#2}\pboxl{0.09\hsize}{$#3$}}&}
	\begin{eqnarray*} 
		\xa{	\partial_t »u - \mathop{\rm div}( M \nabla(\mu \partial_t »u+»w)) \>=\> 	0					}{on}{{\cal T}\times\Omega ,}	\\
		\xa{	»w  \in  \partial \varphi (»u) - \mathop{\rm div}( b_1 (»u,\nabla»u, e ))+ b_2 (»u,\nabla»u, e )						}{on}{{\cal T}\times\Omega ,}	\\
		\xa{	\mathop{\rm div} b_0 (»u,\nabla»u, e ) \>=\>  0,\enspace \enspace \enspace 		 e =\epsilon({\bf u}):=\frac12(D{\bf u}+D{\bf u}^t)	}{on}{{\cal T}\times\Omega ,}	\\[1.5ex]
		\xa{	M \nabla(\mu \partial_t »u+»w)\cdot\vec n	\>=\> 0,\enspace \enspace \enspace    b_1 (»u,\nabla»u, e )\cdot\vec n\>=\>  0	}{on}{{\cal T}\times\partial\Omega  ,}	\\
		\xa{	b_0 (»u,\nabla»u, e )\vec n 			\>=\>   0				}{on}{{\cal T}\times{\Gamma_N},}	\\
		\xa{	{\bf u}  \>=\>   0											}{on}{{\cal T}\times{\Gamma_D},}	\\
		\xa{	»u(0)  \>=\>   u_0 										}{on}{\Omega .}
	\end{eqnarray*}
}%
As a consequence of the mass balance and the boundary conditions, the mean value of the
	concentration does not change over time.
Therefore, after applying a simple shift we can assume $ u_0 $ to have mean value $0$
	which then transfers to all $»u(t)$ for $t\in{\cal T}$.
}

The the remainder of this paper we suppose the following.

\begin{assumption}\label{}
	The domain $\Omega \subset{{\mathbb R}}^N $ is nonempty, open, bounded and connected set with
		Lipschitz boundary $\partial\Omega $.
	The two open subsets ${\Gamma_D},{\Gamma_N}$ of $\partial\Omega $ are disjoint with $\partial\Omega =\overline {{\Gamma_D}}\cup\overline {{\Gamma_N}}$,
		${\Gamma_D}\neq\emp$ and $»G:=\Omega \cup{\Gamma_D}$ is regular in the sense of Gr\"oger~\cite{konni1}.
	${\cal T}=\mathop{]\kern1pt 0,T\kern1pt[},\enspace  T>0$ is a bounded (time) interval and $0< m_1 \le m_2 $ and $ q_0 >2$ are
		real constants.         
\end{assumption}

The following regularity result is due to Gr\"oger~\cite{konni1} and applies to regular sets
	$»G$ being the union of a domain $\Omega $ and a part ${\Gamma_D}$ of its boundary $\partial\Omega $,
	where the latter serves as the Dirichlet boundary part.
Before stating his result, spaces used in the formulation are introduced.

\begin{definition}\label{}
	Let
		$ H :=\{»u\in L^2(\Omega ) »: \int_\Omega  u=0 \}$ with the induced $L^2$-inner product and
		$ V  »:=  H^1(\Omega )\cap H$	with the inner product $ \iKlammerA{u}{v}  \iSub{ V }:= \iKlammerA{\nabla u}{\nabla v}  \iSub{L^2}$.
	We define
		$ U :=\{{\bf u}\in H^1(\Omega ;{{\mathbb R}}^N )»: {\bf u}|_{\Gamma_D}=0\}$,
		where ${\bf u}|_{\Gamma_D}$ is understood in the sense of traces of \,${\bf u}$ on ${\Gamma_D}\subset\partial\Omega $
		and with the induced norm of $H^1(\Omega ;{{\mathbb R}}^N )$.
	The mapping $\epsilon$ is given by
	\[
		\epsilon: U  \rightarrow   L^2(\Omega ;{{\mathbb R}}^{N\times N}),\enspace \enspace  \epsilon({\bf u}):=\displaystyle\f12\,(D{\bf u}+D{\bf u}^t). 
	\]
	We equip the range space $ Y :=\epsilon( U )$ with the norm of $L^2(\Omega ;{{\mathbb R}}^{N\times N})$.
	Moreover, for $1\le p\le\infty$, the space $ W^{1,p}_0(»G;{{\mathbb R}}^M ) $ is defined to be the closure~of
		\[
			\{{\bf u}|_{\Int»G} »: {\bf u}\in C^\infty_c({{\mathbb R}}^N ;{{\mathbb R}}^M ),\enspace  \mathop{\rm supp}{\bf u}\cap(\overline {»G}\sm»G) = \emp \}
		\]
		in the usual Sobolev spaces $W^{1,p}(\Int»G;{{\mathbb R}}^M )$ and
		$ W^{-1,p}(»G;{{\mathbb R}}^M ) :=(W^{1,p'}_0(»G;{{\mathbb R}}^M ))^*$ for the conjugated exponent $p'$ given by
		$\frac1p+\frac1{p'}=1$ (and using the convention $\frac1\infty:=0$).
\end{definition}

Let us agree to simply write $ ||\iBlock{x}||  _{ H } $ for $ ||\iBlock{x}||  _{L^2} $ even if $»x\not\in X $.
Now we are in the position to state Gr\"oger's regularity result adapted to our situation.

\begin{proposition}\label{p4:Konni}
	Let $»b:»G\times{{\mathbb R}}^{N\times N}\rightarrow  {{\mathbb R}}^{N\times N}$ such that 
		\begin{eqnarray*} 
			&	»x\mapsto »b(»x,0)			\,\in\,		L^{ q_0 }(»G;{{\mathbb R}}^{N\times N}),		\enspace \enspace \enspace 
				»x\mapsto »b(»x,»v)			\enspace 			\mbox{is measurable},	&	\\
			&	\bbig(»b(»x,»v)-»b(»x,»w)\bbig)\cdot\bbig(»v-»w\bbig)
										\>\ge\> 			m_1  |»v-»w|^2,			\enspace \enspace 
				|»b(»x,»v)-»b(»x,»w)| 	\>\le\>  		m_2  |»v-»w|. 				&				
		\end{eqnarray*}%
	Corresponding to $»b$, the operator $ A : U \rightarrow   U ^*$ is given by
		\[
			 \iKlammerB{ A {\bf u}}{{\bf v}}  \iSub{ U }:=\int_\Omega »b(»x,\epsilon({\bf u})):\epsilon({\bf v})\>dx.
		\]
	Then there exists a constant $ q_1 $ depending only on $»G, m_1  $ and $ m_2  $ with
		$2< q_1 \le q_0 $ such that $ A $ maps the subspace $ W^{1, q_1 }_0(»G;{{\mathbb R}}^M ) $ of $ U $ onto the space $ W^{-1, q_1 }(»G;{{\mathbb R}}^M ) $.
\end{proposition}

\begin{remark}\label{}
	Note that $ U =W^{1,2}_0(»»G;{{\mathbb R}}^{N\times N})$.
	The original result of Gr\"oger~\cite{konni1} is given for scalar functions
		under conditions analog to those given above.
	To this end, he shows that the duality mapping of $W^{1,2}_0(»G;{{\mathbb R}})$ maps the
		subspace $W^{1,p}_0(»G;{{\mathbb R}})$ onto $W^{-1,p}_0(»G;{{\mathbb R}})$ for some $p>2$ and
		then transfers this property to	nonlinear operators. 
	We further note that all arguments of Gr\"oger can be transferred to the vector-valued case
		where the norm $ ||\iBlock{\>.\> }||  _{ U } $ of $ U $ is replaced by the
		equivalent norm $ ||\iBlock{\epsilon(\>.\> )}||  _{ Y } $ (due to Korn's inequality).
\end{remark}


Throughout this section we further assume the following.

\begin{assumption}\label{}
	Let $ C_P :=\sup\{ ||\iBlock{»u}||  _{ H }  : »u\in V , \> ||\iBlock{»u}||  _{ V } \le1\}$ and
		$ q_1 $ with $2< q_1 \le q_0 $ be given as in Proposition~\ref{p4:Konni} for $G'=»G=\Omega \cup{\Gamma_D}$.
	Furthermore, the constants $ q_2 $ and $ q_3 $ with $2\le q_2 , q_3 \le\infty$ are such that $ V $ \!is 
		continuously embedded
		into $L^{ q_2 }(\Omega )$ and compactly embedded into $L^{ q_3 }(\Omega )$.
	Finally, $ q_4 :=\f{ q_3 }2 (1-\f2{ q_1 })$.

\end{assumption}

\begin{remark}\label{..}
	$ C_P $ is the operator norm of the identity as an operator from $ V $ into $ H $,
		which is finite due to Poincar\'e's inequality.
	Moreover, the Sobolev embedding theorem shows that we can choose
		$ q_2 \ge 2$ arbitraryly in ${{\mathbb R}}$ in case of $N=1,2$ and
		$ q_2 =\f{2N}{N-2}$ if $N\ge3$ together with any
		$ q_3 $ such that $2\le q_3 < q_2 $.
\end{remark}

Next, we consider conditions on the functions $ b_0 , b_1 , b_2 $ and $\varphi $.

\begin{definition}\label{}{%
	\def\xa#1#2{©ruck{1cm}{\hfil{\rm #1}\hfil}{#2}}%
	\def\xb{\vrule width0pt\hskip1cm}%
	\def\ya#1{\iWidth{\iUserA}{#1}#1\global\iUserA=\iUserA}%
	\def\yb#1{\pboxl{\iUserA}{$#1$}}%
	\def\yc#1{\pboxl{3.8cm}{$#1$}}%
	\parskip0pt%
	Within the following conditions all inequalities are assumed to hold for all
		$»t\in T ,»x\in\Omega ,»u, u_1 , u_2 \in{{\mathbb R}},»p, p_1 , p_2 \in{{\mathbb R}}^N , e , e_1 , e_2 \in{{\mathbb R}}^{N\times N}:$\vskip4pt
	{\rucky{13mm}{{\rm(H0)}\hfil}{$ M \in{{\mathbb R}}^{N\times N}$ is symmetric and positive-definite and $\mu \ge0$.
				If $\mu >0$, we set $\mu_0 :=0$, otherwise $\mu_0 :=1$.}}
	\vskip0.8em 
	{\rucky{13mm}{{\rm(H1)}\hfil}{$\varphi :{{\mathbb R}}\rightarrow  \overline {{{\mathbb R}}}$ is a convex, lower-\discretionary{}{}{}semi\discretionary{-}{}{}con\discretionary{-}{}{}tin\discretionary{-}{}{}u\discretionary{-}{}{}ous, proper functional.}}
	\vskip.3em 
	{\rucky{13mm}{{\rm(H1a)}\hfil}{$\varphi \in C^1({{\mathbb R}})$ is convex and $\varphi (r)\le C(r^2+1)$ for all $r\in{{\mathbb R}}$ for some $C>0$.}}
	\vskip0.8em 
	{\rucky{13mm}{{\rm(H2)}\hfil}{$ b_1 : T \times  \Omega \times  {{\mathbb R}}\times  {{\mathbb R}}^N \times {{\mathbb R}}^{N\times N}\rightarrow    {{\mathbb R}}^N $ \>is a Carath\'eodory function with				\\[2pt]
					\xb $\ya{ ( b_1 (»t,»x,»u, p_1 , e )- b_1 (»t,»x,»u, p_2 , e ))\cdot( p_1 - p_2 ) }  \enspace \ge\enspace    \alpha_{ b_1 ,p}| p_1 - p_2 |^2$,					\\[2pt]
					\xb $\yb{ | b_1 (»t,»x,»u,»p, e_1 )- b_1 (»t,»x,»u,»p, e_2 )|               }  \enspace \le\enspace    \beta_{ b_1 ,e}| e_1 - e_2 |$,						\\[2pt]
					\xb $\yc{ | b_1 (»t,»x,»u,»p, e )|^2                      }  \enspace \le\enspace    g (»t,»x) + C_{ b_1 ,u}|»u|^2+ C |»p|^2+ C_{ b_1 ,e}| e |^2$	\\[4pt]
				for some constants $\alpha_{ b_1 ,p}>0,\> \beta_{ b_1 ,e}, C , C_{ b_1 ,u}, C_{ b_1 ,e}\ge0$ and $ g \in {\cal L}^1( {\cal T}\times\Omega  )$.		\\[2pt]
				$ b_2 : T \times  \Omega \times  {{\mathbb R}}\times  {{\mathbb R}}^N \times {{\mathbb R}}^{N\times N}\rightarrow    {{\mathbb R}}$ \>is a Carath\'eodory function with	\\[2pt]
					\xb $\yb{ | b_2 (»t,»x,»u, p_1 , e )- b_2 (»t,»x,»u, p_2 , e ))| }  \enspace \le\enspace    \beta_{ b_2 ,p}| p_1 - p_2 |$,\\[2pt]
					\xb $\yb{ | b_2 (»t,»x,»u,»p, e_1 )- b_2 (»t,»x,»u,»p, e_2 )|  }  \enspace \le\enspace    \beta_{ b_2 ,e}| e_1 - e_2 |$,\\[4pt]
					\xb $\yc{ | b_2 (»t,»x,»u,»p, e )|^2         }  \enspace \le\enspace    g (»t,»x) + C_{ b_2 ,u}|»u|^2+ C |»p|^2+ C_{ b_2 ,e}| e |^2$\\[2pt]
				for some constants $\beta_{ b_2 ,p},\beta_{ b_2 ,e}, C , C_{ b_2 ,u}, C_{ b_2 ,e}\ge0$ and $ g \in {\cal L}^1( {\cal T}\times\Omega  )$.}}
	\vskip0.8em 
	{\rucky{13mm}{{\rm(H3)}\hfil}{$ b_0 : T \times  \Omega \times  {{\mathbb R}}\times  {{\mathbb R}}^N \times {{\mathbb R}}^{N\times N}\rightarrow    {{\mathbb R}}^{N\times N}$ \>is a Carath\'eodory function~with\\[2pt]
					\xb $\yb{ ( b_0 (»t,»x,»u,»p, e_1 )- b_0 (»t,»x,»u,»p, e_2 )):( e_1 - e_2 ) }  \enspace \ge\enspace    \alpha_{ b_0 ,e}| e_1 - e_2 |^2$,\\[2pt]
					\xb $\yb{ | b_0 (»t,»x,»u,»p, e_1 )- b_0 (»t,»x,»u,»p, e_2 ))| }  \enspace \le\enspace    \beta_{ b_0 ,e}| e_1 - e_2 |$,\\[2pt]
					\xb $\yb{ | b_0 (»t,»x,»u, p_1 , e )- b_0 (»t,»x,»u, p_2 , e ))| }  \enspace \le\enspace    \beta_{ b_0 ,p}| p_1 - p_2 |$,\\[2pt]
					\xb $\yc{ | b_0 (»t,»x,»u,»p, e )|^2                     }  \enspace \le\enspace   g (»t,»x)  + C_{ b_0 ,u}|»u|^2+ C |»p|^2+ C_{ b_0 ,e}| e |^2$\\[2pt]
				for some constants $\alpha_{ b_0 ,e}>0,\> \beta_{ b_0 ,e},\beta_{ b_0 ,p}, C , C_{ b_0 ,u}, C_{ b_0 ,e}\ge0$ and $ g \in {\cal L}^1( {\cal T}\times\Omega  )$.\\
				Moreover, the matrix $ b_0 (»t,»x,»u,»p, e )$ is symmetric and continuous in $»t$ uniformly in~$(»x,»u,»p, e )$.}}
	\vskip.3em 
	{\rucky{13mm}{{\rm(H3a)}\hfil}{Condition {\rm(H3)} is satisfied and furthermore it holds\\[2pt]
					\xb $\pboxl{62mm}{$ | b_0 (»t,»x, u_1 ,»p, e )- b_0 (»t,»x, u_2 ,»p, e )|  $}  \enspace \le\enspace    \gamma_{ b_0 ,u}\,(| e |+1)\,| u_1 - u_2 |^{ q_4 }$,\\[4pt]
					\xb $\pboxl{50mm}{$ | b_0 (»t,»x,»u,»p, e )|^{ q_0 }   $}  \enspace \le\enspace    g (»t,»x)  + C (|»u|^{ q_0 }+|»p|^2+| e |^{ q_0 })$,\\[4pt]
				for some $\gamma_{ b_0 ,u}, C \ge0$ and $g\in{\cal L}^1( {\cal T}\times\Omega  )$.
				}}
	\vskip0.8em 
	{\rucky{13mm}{{\rm(H4)}\hfil}{$\alpha_{ b_1 ,p}- C_P \beta_{ b_2 ,p}>0$. Furthermore, we have \,$ m_1 \le \alpha_{ b_0 ,e} \le \beta_{ b_0 ,e} \le m_2 $.}}
	\vskip.3em 
	{\rucky{13mm}{{\rm(H4a)}\hfil}{There exists a constant $ c_a >0$ such that with $\varphi _a $ defined by\\[2pt]
				$\varphi _a :=\bigg[ ( c_a +1) \bigg( C_{ b_1 ,u}+ \f1{\alpha_{ b_0 ,e}} C_P C_{ b_1 ,e} C_{ b_0 ,u}\bigg) »+ (1+\f1{ c_a }) C_P \bigg( C_{ b_2 ,u}+\f1{\alpha_{ b_0 ,e}} C_P C_{ b_2 ,e} C_{ b_0 ,u}) \bigg) \bigg]^{1/2}$    \\[2pt]
				it holds: \enspace \enspace \enspace \enspace  $(\alpha_{ b_1 ,p} - C_P \,\beta_{ b_2 ,p})\alpha_{ b_0 ,e}-(\beta_{ b_1 ,e}+ C_P \,\beta_{ b_2 ,e})\beta_{ b_0 ,p} - \varphi _a   «>0$.}}
	\vskip0.8em 
	{\rucky{13mm}{{\rm(H5)}\hfil}{$ u_0 \in H^1(\Omega )$ with $\int_\Omega u_0 \,dx=0$ and \,$\varphi \circ u_0 \in L^1(\Omega )$.}}
	\vskip15pt 
	{\rucky{13mm}{{\rm(H)}\hfil}{Conditions {\rm(H0)}, {\rm(H1)}, {\rm(H2)}, {\rm(H3)}, {\rm(H4)} and {\rm(H5)} are satisfied.}}
	{\rucky{13mm}{{\rm(Ha)}\hfil}{Conditions {\rm(H0)}, {\rm(H1a)}, {\rm(H2)}, {\rm(H3a)}, {\rm(H4a)} and {\rm(H5)} are satisfied.}}
	\vskip.3em 
}\end{definition}

Here, by Carath\'eodory function we mean a function that is measurable as a function of $(t,x)$ and
	continuous in the other arguments.

\begin{remark}\label{}\nix \\
  \ifnum\the\remarkcc=0\else \\\fi\advance\remarkcc by 1{\bf \the\remarkcc.}\enspace  In condition {\rm(H)} we collect Lipschitz and growth conditions that are needed in order
				to define the operators involved in our weak formulation of problem~{\rm(P)} and
				to apply Theorem~\ref{t1:SemReduct} later on,
				whereas under {\rm(Ha)} we are able to show that
				the weak formulation indeed possesses a solution.\ifnum\the\remarkcc=0\else \\\fi\advance\remarkcc by 1{\bf \the\remarkcc.}\enspace  Under {\rm(H)} the operator $\partial \varphi $ in general is multi-valued.
			 The differentiability assumption in~{\rm(H2a)} could be omitted leading to
				multi-valued pseu\discretionary{-}{}{}domono\discretionary{-}{}{}tone operators later on that can be handled by generalizing
				the results of the preceding »sections.
			 Nevertheless, for the sake of simplicity, we restrict our discussion to
				single-valued pseu\discretionary{-}{}{}domono\discretionary{-}{}{}tone operators here.
				$ m_1 :=\alpha_{ b_0 ,e}$ and $ m_2 :=\beta_{ b_0 ,e}$, if then~{\rm(H3a)} is satisfied.
			 Note that $ q_1 $ and hence $ q_4 $ depend on $ m_1 $ and $ m_2 $.
				(Therefore we fix all these constants in advance.)\ifnum\the\remarkcc=0\else \\\fi\advance\remarkcc by 1{\bf \the\remarkcc.}\enspace  The mappings $ b_1 $ and $ b_2 $ together as well as $ b_0 $ give rise to
				operators $\widetilde  A$ and $\widetilde  B$, respectively, which correspond to the operators
				given in Section~\ref{s2:}.
             Condition~{\rm(H4)} is the strong monotonicity of $\widetilde  B$ and {\rm(H4a)} is a tightening of~(A3.4).\end{remark}

In order to define a (appropriate weak) solution to problem~{\rm(P)} in the framework of Section~\ref{s3:}
	consider the following operators.

\begin{definition}\label{}
	We define $»F\in L ( V ; V ^*),\enspace   \iKlammerB{»F»u}{»v}  \iSub{ V } »:=  \iKlammerA{ M \nabla»u}{\nabla»v} _{L^2(\Omega ;{{\mathbb R}}^N )}$ for $»u,»v\in V $,
	\[
			I_H 	:=	{\rm Id}_{ V \rightarrow   H }									,\enspace \enspace 
			»I	:=	I_H ^* J_H I_H  \in L ( V ; V ^*)						,\enspace \enspace 
			E _1 	:=	\mu {\rm Id}_{ H } + I_H »F^{-1} I_H ^* J_H 		\>\in\>  L ( H ; H )	.
	\]
	Let $ E _2 \in L( H )$ be the (positive, symmetric) root of $ E _1 $ and
	\[
			K 	:=	E _2 I_H 							\>\in\>  L ( V ; H )	,\enspace \enspace \enspace 
			E 	:=	K ^* J_H K 						\>\in\>  L ( V ; V ^*)	.
	\]
	Corresponding to these spaces and operators let ${\cal  V },{\cal  H },{\cal W},{\cal E},{\cal  K },{\cal  L }$ and ${\cal  L }_{»»h}$
		be given as in Section~\ref{s3:} and define ${\cal U}:=L^2( T ; U )$.
\end{definition}

\begin{remark}\label{r4:E}
	The operator $»F$ corresponds to the mapping $-\mathop{\rm div}( M \nabla\>.\> )$ with natural boundary conditions
		and is positive-definite and symmetric.
	These properties of $»F$ also transfer to $ E _1 $ and $ E $ and it holds
	\[
		E   \>=\>   I_H ^* E _2 ^* »JK E _2 I_H   \>=\>   \mu »I + »I»F^{-1}»I.
	\]
	Note that $ K $ has dense range.
	Indeed, in the case of $\mu =0$ the operator $ E _1 $ is the composition of operators with dense range.
	If $\mu >0$ then $ E _1 $ is even surjective since it is monotone, continuous and coercive.
	Therefore, also $ E _2 $ and $ K $ have dense range in $ H $.	
\end{remark}

As already done, we identify the dual of the space $L^2({\cal T};X)$ for a reflexive Banach space $ X $
	with the space $L^2({\cal T};X^*)$.
Now, we introduce operators related to the functions $ b_0 , b_1 , b_2 $ and $\varphi $.
The Carath\'eodory property and the growth conditions ensure that these operators are indeed mappings between the given spaces.

\begin{definition}\label{..}
	Suppose that conditions {\rm(H1)}, {\rm(H2)} and {\rm(H3)} hold. We set
		{\def\xa#1#2#3#4{&\pboxl{6mm}{$#1$} : \pboxl{38mm}{$#2$}\enspace \enspace \enspace  \pboxl{9.5cm}{\pboxl{3.9cm}{$#3$} «:= \pboxl{5.0cm}{$#4$}}&}%
		\def\xb#1#2#3{  &\pboxl{6mm}{$#1$} : \pboxl{38mm}{$#2$}\enspace \enspace \enspace  \pboxl{9.5cm}{\pboxl{3.9cm}{$#3$} \pboxc{11mm}{} \pboxl{5.0cm}{}}&}%
		\begin{eqnarray*} 
			\xa{ \xhh{\widetilde  B^Y}{B} }{  {\cal T}\times V \times V \times Y   \rightarrow   Y ^*, }{  \iKlammerB{\xhh{\widetilde  B^Y}{B}(»t, u_1 , u_2 , e )}{ e _1 }  \iSub{ Y } }{ \int_\Omega  b_0 (t,»x, u_2 ,\nabla u_1 , e ): e _1  \> dx,       }	\\[0.5ex]
			\xa{ \xhh{\widetilde  B^1}{B} }{  {\cal T}\times V \times V \times Y   \rightarrow   V ^*, }{  \iKlammerB{\xhh{\widetilde  B^1}{B}(»t, u_1 , u_2 , e )}{»v}  \iSub{ V }  }{ \int_\Omega  b_1 (t,»x, u_2 ,\nabla u_1 , e )\cdot\nabla»v \> dx, }	\\[0.5ex]
			\xa{ \xhh{\widetilde  B^2}{B} }{  {\cal T}\times V \times V \times Y   \rightarrow   V ^*, }{  \iKlammerB{\xhh{\widetilde  B^2}{B}(»t, u_1 , u_2 , e )}{»v}  \iSub{ V }  }{ \int_\Omega  b_2 (t,»x, u_2 ,\nabla u_1 , e )\,»v \> dx,       }	\\[0.5ex]
			\xb{ \xhh{\widetilde  B^X}{B} }{  {\cal T}\times V \times V \times Y   \rightarrow   V ^*, }{ \xhh{\widetilde  B^X}{B}(»t, u_1 , u_2 , e ) «:=  \xhh{\widetilde  B^1}{B}(»t, u_1 , u_2 , e )+\xhh{\widetilde  B^2}{B}(»t, u_1 , u_2 , e ),      }    
		\end{eqnarray*}}%
	 together with
		{\def\xa#1#2#3#4{&\pboxl{5mm}{$#1$} : \pboxl{3.2cm}{$#2$}\enspace \enspace \enspace \enspace  \pboxl{1.9cm}{$#3$} «:= \pboxl{2.3cm}{$#4$}&}%
		\begin{eqnarray*} 
			\xa{ \xhh{B^X}{B} }{ {\cal T}\times V \times Y \rightarrow   V ^*, }{ \xhh{B^X}{B}(»t,»u, e )  }{ \xhh{\widetilde  B^X}{B}(»t,»u,»u, e ),	}    \\[0.5ex]
			\xa{ \xhh{B^Y}{B} }{ {\cal T}\times V \times Y \rightarrow   Y ^*, }{ \xhh{B^Y}{B}(»t,»u, e )  }{ \xhh{\widetilde  B^Y}{B}(»t,»u,»u, e ).	}
		\end{eqnarray*}}%
	Moreover, the operators $\xhh{B^X}{B}$ and $\xhh{B^Y}{B}$ will be extended by
		{\def\xa#1#2#3#4{&\pboxl{5mm}{$#1$} : \pboxl{2.35cm}{$#2$}\enspace \enspace \enspace \enspace  \pboxl{2.7cm}{$#3$} «:= \pboxl{4.5cm}{$#4$}&}%
		\begin{eqnarray*} 
			\xa{ \xhh{{\cal B}^X}{B}  }{ {\cal  V }\times{\cal Y}\rightarrow  {\cal  V }^*, }{  \iKlammerB{\xhh{{\cal B}^X}{B}(»u, e )}{»v}  \iSub{{\cal  V }}  }{ \int_{\cal T}   \iKlammerB{\xhh{B^X}{B}(t,»u, e )}{»v}  \iSub{ V } \> dt,  }    \\[0.5ex]
			\xa{ \xhh{{\cal B}^Y}{B}  }{ {\cal  V }\times{\cal Y}\rightarrow  {\cal Y}^*, }{  \iKlammerB{\xhh{{\cal B}^Y}{B}(»u, e )}{ e _1 }  \iSub{{\cal Y}} }{ \int_{\cal T}   \iKlammerB{\xhh{B^Y}{B}(t,»u, e )}{ e _1 }  \iSub{ Y } \> dt. }
		\end{eqnarray*}}%
		to operators on ${\cal  V }\times{\cal Y}$.
	Note that again the dependence of \>$»u,\nabla»u,»v, e $ and $ e _1 $ on $»x\in\Omega $ and $»t\in T $
		was suppressed in this notation.
	{\def\xa#1#2#3{\pboxc{35mm}{$#1$}\pboxl{10mm}{$#2$}:=\enspace \pboxl{6cm}{$#3$}}%
	Moreover, we define the functionals
		\begin{eqnarray*} 
				\xa{  Q : V \rightarrow  \overline {{{\mathbb R}}},   }{ Q (»u)  }{
					\begin{cases}
						\textstyle\int\limits_\Omega { \varphi  \circ»u  }\enspace \enspace 	& \text{if $\varphi  \circ»u\in L^1(\Omega )$},\\
						+\infty									& \text{otherwise.}
					\end{cases}						}\\[4pt]
 				\xa{  {\cal Q}:{\cal  V }\rightarrow  \overline {{{\mathbb R}}},  }{ {\cal Q}(»u)  }{
					\begin{cases}
						\textstyle\int\limits_{\cal T}{ Q \circ»u  }\enspace \enspace 	& \text{if $ Q \circ»u\in L^1({\cal T})$},\\
						+\infty									& \text{otherwise.}
					\end{cases}						}
		\end{eqnarray*}%
		and the operator ${\cal A} »:= \partial{\cal Q} \>\subset\> {\cal  V }\times{\cal  V }^*$.
	}
\end{definition}

Now we are in the position to introduce our concept of weak solutions for problem~{\rm(P)}.

\begin{definition}[{\rm Weak formulation}]\label{d4:Weak}
	A tuple $(»u,{\bf u})\in{\cal W}\times{\cal U}$ is called a weak solution to problem~{\rm(P)} if
		for $ e :=\epsilon({\bf u})\in Y $ it holds
	\[
		{\cal  L }»u+{\cal A}»u+\xhh{{\cal B}^X}{B}(»u, e ) \>\ni\> 0,\enspace \enspace \enspace  \xhh{{\cal B}^Y}{B}(»u, e )»=0,\enspace \enspace \enspace  ({\cal  K }»u)(0)»= K u_0 .
	\]
\end{definition}

\begin{remark}\label{..}\nix\\
	\ifnum\the\remarkcc=0\else \\\fi\advance\remarkcc by 1{\bf \the\remarkcc.}\enspace  %
		By virtue of Proposition~\ref{p3:PartInt}, the images of functions of ${\cal W}$ under the mapping
		${\cal  K }$ can be regarded as elements of $C(\overline {{\cal T}}; H )$.
		Therefore, $({\cal  K }»u)(0)\in H $ is well defined and the condition $({\cal  K }»u)(0) »= K u_0 $ meaningful.\ifnum\the\remarkcc=0\else \\\fi\advance\remarkcc by 1{\bf \the\remarkcc.}\enspace  %
		It is not hard to show that if $»u,»w$ and ${\bf u}$ are sufficiently smooth in the sense of
		Sobolev spaces, they are strong or even classical solutions to problem~{\rm(P)}.\ifnum\the\remarkcc=0\else \\\fi\advance\remarkcc by 1{\bf \the\remarkcc.}\enspace  %
		In our weak formulation of problem~{\rm(P)} we only require ${\cal E}»u$ to have generalized derivatives
			within ${\cal  V }^*$, not $»I»u$ itself.
		This relaxation of the regularity requirements together with the linearity
			of $»u'-\mathop{\rm div}( M \nabla»w)=0$ allow us to treat problem~{\rm(P)} with the techniques of Section~\ref{s3:}.
		Note that ${\cal A}»u+\xhh{{\cal B}^X}{B}(»u, e )$ only contributes space derivatives up to second order.
		The remaining ones are 'hidden' in the operator $ E $.
		Roughly speaking, the chemical potential $»w$ only attains values in $ V ^*$,
			but values in $ V $ are needed in order to use the standard weak formulation of
			the diffusion equation $\partial_t »u »- \mathop{\rm div}( M \nabla(\mu \partial_t »u+»w))=0$.
		Therefore, we apply the operator $»I»F^{-1}$ to this equation, eliminate $»w$ and
			use the resulting equation as a new weak formulation.\end{remark}

In order to show that problem~{\rm(P)} possesses a weak solution we show 
	firstly, that	for $»t\in{\cal T}$ the equation $\xhh{B^Y}{B}(»t,»u, e )=0$ has a unique solution $ e = e (»u)\in Y $
	for every $»u\in V $
	and secondly, that mapping $»u\mapsto \xhh{B^X}{B}(»u, e (»u))$ is pseu\discretionary{-}{}{}domono\discretionary{-}{}{}tone.
Consequently, Theorem~\ref{t2:Exist} will guarantee the existence of weak solutions.

\begin{lemma}\label{lab.lemma.red.lip}{\def\xa#1#2#3{\pboxl{8cm}{$#1$}\pboxc{1cm}{$#2$}\pboxl{5.1cm}{$#3$}}%
	\def\xb#1#2#3{\pboxl{8cm}{$#1$}\pboxc{1cm}{$#2$}\pboxl{5.1cm}{$#3$}}%
	Let {\rm(H2)} and {\rm(H3)} be satisfied. Then it follows that
	\begin{eqnarray*} 
			\xa{  \iKlammerB{\xhh{\widetilde  B^X_t}{B}( u_1 ,»u, e )-\xhh{\widetilde  B^X_t}{B}( u_2 ,»u, e )}{ u_1 - u_2 }  \iSub{ V }     }{ \ge }{ (\alpha_{ b_1 ,p}- C_P \beta_{ b_2 ,p}) ||\iBlock{ u_1 - u_2 }||  _{ V } ^2,    }\\[0.5ex]
			\xa{  \iKlammerB{\xhh{\widetilde  B^Y_t}{B}( u_1 , u_2 , e _1 )-\xhh{\widetilde  B^Y_t}{B}( u_1 , u_2 , e _2 )}{ e _1 - e _2 }  \iSub{ Y } }{\ge }{ \alpha_{ b_0 ,e} ||\iBlock{ e _1 - e _2 }||  _{ Y } ^2,               }\\[7pt]
			\xb{  ||\iBlock{\xhh{\widetilde  B^X_t}{B}( u_1 , u_2 , e _1 )-\xhh{\widetilde  B^X_t}{B}( u_1 , u_2 , e _2 )}||  _{ V ^*}          }{ \le }{ (\beta_{ b_1 ,e}+ C_P \beta_{ b_2 ,e})  ||\iBlock{ e _1 - e _2 }||  _{ Y } ,     }\\[0.5ex]
			\xb{  ||\iBlock{\xhh{\widetilde  B^Y_t}{B}( u_1 , u_2 , e _1 )-\xhh{\widetilde  B^Y_t}{B}( u_1 , u_2 , e _2 )}||  _{ V ^*}          }{ \le }{ \beta_{ b_0 ,e}  ||\iBlock{ e _1 - e _2 }||  _{ Y } ,                }\\[0.5ex]
			\xb{  ||\iBlock{\xhh{\widetilde  B^Y_t}{B}( u_1 ,»u, e )-\xhh{\widetilde  B^Y_t}{B}( u_2 ,»u, e )}||  _{ V ^*}              }{ \le }{ \beta_{ b_0 ,p}  ||\iBlock{ u_1 - u_2 }||  _{ V }                  }
	\end{eqnarray*}%
	for all $»t\in T ,»u, u_1 , u_2 \in V $ \!and $ e , e _1 , e _2 \in Y $.
	In case of {\rm(H3b)} we also have
	\begin{eqnarray*} 
			\xb{  ||\iBlock{\xhh{\widetilde  B^Y_t}{B}(»u, u_1 , e )-\xhh{\widetilde  B^Y_t}{B}(»u, u_2 , e )}||  _{ Y ^*}              }{ \le }{ C_P \beta_{ b_0 ,u}  ||\iBlock{ u_1 - u_2 }||  _{ V }               }
	\end{eqnarray*}%
}\end{lemma}

{\em Proof}. 
	Exemplarily, we show the strong monotonicity of $\xhh{\widetilde  B^X_t}{B}$ in the first argument and the Lipschitz continuity in the last component.
	The other inequalities can be proven similarly. 
	To this end, suppose that $»u, u_1 , u_2 \in V $ and $ e \in Y $.
	Due to the definition of $ C_P $ we have
	\[
			 ||\iBlock{ u_1 - u_2 }||  _{ H }   \>\le\>    C_P \, ||\iBlock{\nabla( u_1 - u_2 )}||  _{ H }  »= C_P \, ||\iBlock{ u_1 - u_2 }||  _{ V } .
	\]
	The Cauchy-Schwarz inequality and {\rm(H2)} yield
	\begin{eqnarray*} 
			&     &   \hskip-2cm   \iKlammerB{\xhh{\widetilde  B^X_t}{B}( u_1 ,»u, e )-\xhh{\widetilde  B^X_t}{B}( u_2 ,»u, e )}{ u_1 - u_2 }  \iSub{ V }\\[0.5ex]
			& »=  &  \int_\Omega    \bbig( b_1 (»t,»x,»u,\nabla u_1 , e )- b_1 (»t,»x,»u,\nabla u_2 , e ) \bbig)\cdot\nabla( u_1 - u_2 )   \> dx  \\
			&     &  +\> \int_\Omega    \bbig( b_2 (»t,»x,»u,\nabla u_1 , e )- b_2 (»t,»x,»u,\nabla u_2 , e ) \bbig)\cdot( u_1 - u_2 )  \> dx  \\[0.5ex]
			& \>\ge\>   &  \alpha_{ b_1 ,p}\, ||\iBlock{\nabla( u_1 - u_2 )}||  _{ H } ^2  »-  \beta_{ b_2 ,p}\, ||\iBlock{\nabla( u_1 - u_2 )}||  _{ H }   ||\iBlock{ u_1 - u_2 }||  _{ H }   \\[0.5ex]
			& \>\ge\>   &  (\alpha_{ b_1 ,p}- C_P \,\beta_{ b_2 ,p})\, ||\iBlock{ u_1 - u_2 }||  _{ V } ^2.
	\end{eqnarray*}%
	In order to show the Lipschitz continuity of $\xhh{\widetilde  B^X_t}{B}$ in the last argument we estimate
	\begin{eqnarray*} 
			&     &   \hskip-2cm   \iKlammerB{ \xhh{\widetilde  B^X_t}{B}( u_1 , u_2 , e _1 )-\xhh{\widetilde  B^X_t}{B}( u_1 , u_2 , e _2 ) }{ »u }  \iSub{ V }\\[0.5ex]
			& »=  &  \int_\Omega    \bbig( b_1 (»t,»x, u_2 ,\nabla u_1 , e _1 )- b_1 (»t,»x, u_2 ,\nabla u_1 , e _2 ) \bbig)\cdot\nabla»u   \> dx  \\
			&     &  +\> \int_\Omega    \bbig( b_2 (»t,»x, u_2 ,\nabla u_1 , e _1 )- b_2 (»t,»x, u_2 ,\nabla u_1 , e _2 ) \bbig)»u\,  \> dx  \\[0.5ex]
			& \>\le\>   &  \beta_{ b_1 ,e}\, ||\iBlock{ e _1 - e _2 }||  _{ Y }  ||\iBlock{\nabla»u}||  _{ H }   »+  \beta_{ b_2 ,e}\, ||\iBlock{ e _1 - e _2 }||  _{ Y }   ||\iBlock{»u}||  _{ H }   \\[0.5ex]
			& \>\le\>   &  (\beta_{ b_1 ,e}+ C_P \,\beta_{ b_2 ,e})\, ||\iBlock{ e _1 - e _2 }||  _{ Y }  ||\iBlock{»u}||  _{ H } 
	\end{eqnarray*}%
	for arbitrary $»u, u_1 , u_2 \in V $ and $ e _1 , e _2 \in Y $. Since
	\[
					 ||\iBlock{ \xhh{\widetilde  B^X_t}{B}( u_1 , u_2 , e _1 )-\xhh{\widetilde  B^X_t}{B}( u_1 , u_2 , e _2 ) }||  _{ V ^*} 
		=	\osD{\displaystyle\sup_{»u\in V ,\atop  ||\iBlock{»u}||  _{ V } \le1}}{}   \iKlammerB{ \xhh{\widetilde  B^X_t}{B}( u_1 , u_2 , e _1 )-\xhh{\widetilde  B^X_t}{B}( u_1 , u_2 , e _2 ) }{ »u }  \iSub{ V },	
	\]
	we obtain the desired inequality.
\endproof


\begin{corollary}\label{c4:Ass}
	Suppose {\rm(H)} to be satisfied.
	Then $\xhh{\widetilde  B^X_t}{B}$ and $»\xhh{\widetilde  B^Y_t}{B}$ (as $\widetilde { A }$ and $\widetilde { B }$) satisfy~{\rm(A1)} and {\rm(A2)} of Section~\ref{s2:}
		for every $y_0^*\in Y ^*$.
	Moreover, the constants can be chosen by
	\[
			\alpha_A  := \alpha_{ b_1 ,p} - C_P \,\beta_{ b_2 ,p},		\enspace \enspace \enspace \enspace 
			\beta_A  := \beta_{ b_1 ,e}+ C_P \,\beta_{ b_2 ,e}, 		\enspace \enspace \enspace \enspace 
			\alpha_B  := \alpha_{ b_0 ,e},					\enspace \enspace \enspace \enspace 
			\beta_B  := \beta_{ b_0 ,p}.
	\]
\end{corollary}

As done in Section~\ref{s2:}, we introduce the operators $\xwt{ R }$ and $\xwt{ S }$ which now also depend on~$»t\in T $.

\begin{definition}\label{}{\def\xa#1#2#3#4{\pboxl{6mm}{$#1$}\pboxl{4cm}{$#2$}\pboxl{2.2cm}{$#3$}\pboxl{45mm}{$#4$}}%
	Assume {\rm(H)} to be satisfied.
	Then for every $»t\in T $ we define the operators $\xwt{ R }_t $ and $\xwt{ S }_t $ according to Definition~\ref{d2:RS} and $ y_0^*:=0$ as
	\begin{eqnarray*} 
			\xa{ \xwt{ R }: }{ T \times V \times V \rightarrow   Y ,  }{ \xwt{ R }(»t, u_1 , u_2 ) }{ :=\> (\xhh{\widetilde  B^Y}{B}\hskip-2mm_{t, u_1 , u_2 })^{-1}(0),  } \\[0.5ex]
			\xa{ \xwt{ S }: }{ T \times V \times V \rightarrow   V ^*, }{ \xwt{ S }(»t, u_1 , u_2 ) }{ :=\> (\xhh{\widetilde  B^X_t}{B}( u_1 , u_2 ,\xwt{ R }_t ( u_1 , u_2 )).   }
	\end{eqnarray*}%
	Moreover, let $ B (t):= S _t : V \times V \rightarrow   V ^*$ and ${\cal B}$ the superposition operator
		(Nemytskii operator) of $ B $ given by $({\cal B}»u)(t):= B (t)»u(t)$.
}\end{definition}

\begin{lemma}\label{p4:Growth}{\def\xa#1#2{\parskip4pt\rucky{20mm}{\hskip10mm{\rm #1}\hss}{#2}}%
  Let {\rm(H)} be fulfilled. Then there exist $C>0$ and $»h\in L^1( T )$ such that the following statements
	are fulfilled for all $»t\in T $ and $»u, u_1 , u_2 \in V $:
    \xa{1.}{the mappings $»t\mapsto  R _t »u$ and $»t\mapsto  B _t »u$ are continuous (and hence measurable),}
    \xa{2.}{$   ||\iBlock{\xwt{ R }_t ( u_1 , u_2 )}||  _{ Y } ^2		\enspace \le\enspace     »h(»t)  »+  C  ||\iBlock{ u_1 }||  _{ V } ^2  »+ \f1{\alpha_B }\, C_P C_{ b_0 ,u}  ||\iBlock{ u_2 }||  _{ V } ^2$,}
    \xa{3.}{$   ||\iBlock{\xwt{ S }_t ( u_1 , u_2 )}||  _{ V ^*} ^2	\enspace \le\enspace     »h(»t)  »+  C  ||\iBlock{ u_1 }||  _{ V } ^2  »+ \varphi _a ^2 ||\iBlock{ u_2 }||  _{ V } ^2   $,}
    \xa{4.}{${\cal B}$ is a bounded mapping from ${\cal  V }$ into ${\cal  V }^*$}
}\end{lemma}

{\em Proof}. 
	1. Let $»u\in V , t_0 \in T $ and $ \Ge >0$ be given. We define $ e (»t):= R _t (»u)$.
	The continuity of $ b_0 $ in $t$ implies that $\xhh{\widetilde  B^Y_t}{B}(»u, e ( t_0 ))$ is continuous in $»t$.
	Hence, there exists a $ \delta >0$ such that
		$ ||\iBlock{ \xhh{\widetilde  B^Y_t}{B}(»u, e ( t_0 ))-\xhh{B^Y_{t_0}}{B}(»u, e ( t_0 )) }||  _{ Y ^*}  =  ||\iBlock{ \xhh{\widetilde  B^Y_t}{B}(»u, e ( t_0 )) }||  _{ Y ^*} 
					< \alpha_B  \Ge $
		for all $»t\in T $ with $|»t- t_0 |< \delta $.
	The strong monotonicity of $\xhh{B^Y_{»t,»u}}{B}$ implies the Lipschitz continuity of $(\xhh{B^Y_{»t,»u}}{B})^{-1}$.
	Hence, it holds 
		\[
										 ||\iBlock{ e (»t)- e ( t_0 ) }||  _{ Y } 
					\le	  	\f{1}{\alpha_B }  ||\iBlock{ \xhh{\widetilde  B^Y_t}{B}(»u, e (t))-\xhh{\widetilde  B^Y_t}{B}(»u, e ( t_0 )) }||  _{ Y ^*} 
					=		\f{1}{\alpha_B }  ||\iBlock{ \xhh{\widetilde  B^Y_t}{B}(»u, e ( t_0 )) }||  _{ Y ^*}  < \Ge 
		\]
		for all $»t\in T $ with $|»t- t_0 |< \delta $.
	This proves the continuity of $»t\mapsto  R _t »u$ and hence its measurability.
	Since $\xhh{B^X_t}{B}(»u, e )$ satisfies the Carath\'eodory condition,
	the mapping $»t\mapsto \xhh{B^X_t}{B}(»u,\xwt{ R }_t »u)= B _{»t}»u$ is measurable.

	2.+3. Let $»t\in T ,\, »z=( u_1 , u_2 )\in V \times V $ be given and denote $»y:=\xwt{ R }_t »z$.
	From the strong monotonicity of $\xhh{\widetilde  B^Y_t}{B}$ it follows
		\[
										 ||\iBlock{»y}||  _{ Y } ^2
					\enspace \le\enspace  		\f1{\alpha_B }\,  \iKlammerB{ \osH{\xhh{\widetilde  B^Y_t}{B}(»z,»y)-\xhh{\widetilde  B^Y_t}{B}(»z,»0)}{\vrule width0pt height9pt} }{»y-0}  \iSub{ Y }
					\>\le\>  	\f1{\alpha_B }\, ||\iBlock{\xhh{\widetilde  B^Y_t}{B}(»z,0)}||  _{ Y ^*}  ||\iBlock{»y}||  _{ Y } 
		\]
		because of $\xhh{\widetilde  B^Y_t}{B}(»z,»y)= y_0^*=0$.
	By the growth condition on $ b_0 $, we obtain for some $»h\in L^1( T )$
		\begin{eqnarray*} 
			\alpha_B \,  ||\iBlock{\xwt{ R }_t »z}||  _{ Y }    & \le &   ||\iBlock{\xhh{\widetilde  B^Y_t}{B}(»z,»0)}||  _{ Y ^*}   
						«=    \sup\bigg\{ \iKlammerB{ \osH{\xhh{\widetilde  B^Y_t}{B}(»z,0)}{\vrule width0pt height9pt} }{ e '}  \iSub{ Y }»: e '\in Y ,\enspace   ||\iBlock{ e '}||  _{ Y } \le1  \bigg\} \\
								& \le &   ||\iBlock{ b_0 (»t,»., u_2 ,\nabla u_1 ,0) }||  _{ H }    \\
								& \le &  \bigg( \int_\Omega  »g(»t,»x) \> dx  + C  ||\iBlock{\nabla u_1 }||  _{ H } ^2 + C_{ b_0 ,u}  ||\iBlock{ u_2 }||  _{ H } ^2  \bigg)^{1/2} \\
								& \le &  \bigg( »h(»t)  + C  ||\iBlock{ u_1 }||  _{ V } ^2 + C_P C_{ b_0 ,u}  ||\iBlock{ u_2 }||  _{ V } ^2  \bigg)^{1/2}\hskip-3pt.
		\end{eqnarray*}%
	With this inequality and the growth condition on $ b_1 $ and $ b_2 $ we can similarly estimate
		\begin{eqnarray*} 
			 ||\iBlock{\xwt{ S }_t »z}||  _{ V ^*} 																
							&  =  &   ||\iBlock{\xhh{\widetilde  B^X_t}{B}(»z,\xwt{ R }_t »z)}||  _{ V ^*}   \\[0.5ex]
							& \le &   ||\iBlock{ b_1 (»t,»., u_2 ,\nabla u_1 ,\xwt{ R }_t »z)}||  _{ H }  »+ C_P   ||\iBlock{ b_2 (»t,»., u_2 ,\nabla u_1 ,\xwt{ R }_t »z)}||  _{ H }  \\[0.5ex]
							& \le &  \bigg( »h(»t) »+ C ||\iBlock{\nabla u_1 }||  _{ H } ^2  »+  C_{ b_1 ,u} ||\iBlock{ u_2 }||  _{ H } ^2  »+ C_{ b_1 ,e}  ||\iBlock{\xwt{ R }_t »z}||  _{ H } ^2 \bigg)^{1/2}  \\[-2pt]
							&     & »+ C_P  \bigg( »h(»t) »+ C ||\iBlock{\nabla u_1 }||  _{ H } ^2  »+  C_{ b_2 ,u} ||\iBlock{ u_2 }||  _{ H } ^2  »+ C_{ b_2 ,e}  ||\iBlock{\xwt{ R }_t »z}||  _{ H } ^2 \bigg)^{1/2}  \\
							& \le &   \bigg( »h(»t)   »+  C ||\iBlock{ u_1 }||  _{ V } ^2 »+ \varphi _a ^2 ||\iBlock{ u_2 }||  _{ V } ^2 \bigg)^{1/2}\hskip-3pt.
	\end{eqnarray*}%
	\iHeight\iDa{$\f1c$}\iDepth\iDb{$\f1c$}\def\xa{\vrule width0pt height\iDa depth\iDb}%
	\def\xa{}%
	In the last line the inequality \osD{$\sqrt{\xa a\,}+\sqrt{\xa b\,} \le \sqrt{( c_a +1)a+(1+\f1{ c_a })b\,}$}{}\ for $a,b\ge0$ was used.

	4. The mapping $ B : T \times V \rightarrow   V ^*$ is measurable in $t$ and demicontinuous in $»v$.
	Hence, ${\cal B}»u$ is measurable for every $»u\in{\cal  V }$.
	Moreover, the growth condition of step~2 guarantees that ${\cal B}$ is a bounded operator from ${\cal  V }$ into ${\cal  V }^*$.
\endproof


With the help of this lemma and the bijectivity of $\xhh{B^Y_{»t,»u}}{B}$
	the task of finding a weak solution to problem~{\rm(P)} can be reformulated in the following way.

\begin{corollary}\label{c4:Iff}
	A pair $(»u,{\bf u})\in{\cal W}\times{\cal U}$ is a weak solution to problem~{\rm(P)} if and only if
		$ e (t):= R (»t,»u(»t))$ satisfies ${\bf u}(t)=\epsilon^{-1}( e (t))$ and
		$»u\in{\cal W}$ is a solution to
		\[
			({\cal  L }+{\cal A}+{\cal B})»u \>\ni\> 0,\enspace \enspace   ({\cal  K }»u)(0) »= K u_0 .
		\]
\end{corollary}

\begin{lemma}\label{p4:Ssc}
	Assume {\rm(H2)},\,{\rm(H3a)} and {\rm(H4)} to be fulfilled. Then
		$
			\xhh{\widetilde  B^X_t}{B}: V \times V _\omega \times Y  \>\rightarrow  \> V ^*
		$
		is continuous for all $»t\in T $.
	Furthermore, for $»u_n \mathop{\rightharpoonup}»u$ in $ V $ it holds
	\[
			»\xhh{\widetilde  B^Y_t}{B}(»v,»u_n , e ) \mathop{\relbar\joinrell\rightarrow} »\xhh{\widetilde  B^Y_t}{B}(»v,»u, e ) \mbox{\enspace \enspace  in $ Y ^*$}
	\]
	for all $»t\in T ,\> »v,»v_1 ,»v_2 \in V $ and every solution $ e \in Y $ to $»\xhh{\widetilde  B^Y_t}{B}(»v_1 ,»v_2 , e )=0$.
\end{lemma}

{\em Proof}. 
	The continuity of $\xhh{\widetilde  B^X_t}{B}$ is a direct consequence of the growth conditions on $ b_1 $ and $ b_2 $ and
		the compact embedding of $ V _\omega $ into $ H $.
	Assume that $»v,»v_1 ,»v_2 \in V $ are given, $»u_n \mathop{\rightharpoonup}»u$ in $ V $ and that $ e $
		is a solution to $»\xhh{\widetilde  B^Y_t}{B}(»v_1 ,»v_2 , e )=0$.
	By {\rm(H2)}, the mapping $ e '\mapsto  b_0 (»t,»x,»v_2 (»x),\nabla»v_1 (»x), e ')$ is strongly monotone and
		Lipschitz continuous from ${{\mathbb R}}^{N\times N}$ into itself independently of $(»t,»x)\in T \times\Omega $.
	Furthermore, due to {\rm(H3a)}, $»x\mapsto  b_0 (»t,»x,»v_2 (»x),\nabla»v_1 (»x),0)\in L^{ q_0 }(\Omega ;{{\mathbb R}}^{N\times N})$ for all $»t\in T $.
	Consequently, Proposition~\ref{p4:Konni} implies $ e \in L^{ q_1 }(\Omega ;{{\mathbb R}}^{N\times N})$ for every $»t\in T $.
	Moreover, the convergence $»u_n \mathop{\rightharpoonup}»u$ in $ V $ yields $»u_n \mathop{\rightarrow }»u$ in $L^{ q_3 }(\Omega )$.
	Therefore, by {\rm(H3a)} and H\"older's inequality we get for all $»t\in T $ and $ e '\in Y $
	\begin{eqnarray*} 
			&& \hskip-2cm   \iKlammerB{ »\xhh{\widetilde  B^Y_t}{B}(»v,»u_n , e )-»\xhh{\widetilde  B^Y_t}{B}(»v,»u, e ) }{ e ' }  \iSub{ Y }                            	\\[0.5ex]
			&  =  & 	\int_\Omega   \big[ b_0 (»t,»x,»u_n ,»v, e )- b_0 (»t,»x,»u,»v, e ) \big]: e ' \> dx   		\\[0.5ex]
			& \le & 	 ||\iBlock{ e '}||  _{ Y }  \int_\Omega   \bbig| b_0 (»t,»x,»u_n ,»v, e )- b_0 (»t,»x,»u,»v, e ) \bbig|^2 \> dx	\\[0.5ex]
			& \le & 	\gamma_{ b_0 ,u}^2  ||\iBlock{ e '}||  _{ Y }   \int_\Omega   | u_1 - u_2 |^{2 q_4 }(| e |+1)^2 \> dx                   	\\[0.5ex]
			& \le & 	C\,      ||\iBlock{ e '}||  _{ Y }    ||\iBlock{ u_1 - u_2 }||  _{{L^{ q_3 }(\Omega )}} ^{2 q_4 }  ( ||\iBlock{ e }||  _{L^{ q_1 }(\Omega )} ^2+1)
	\end{eqnarray*}%
	since
		$(\f{ q_3 }{2 q_4 })^{-1} + (\f{ q_1 }{2})^{-1} = \f{ q_1 -2}{ q_1 } + \f{2}{ q_1 } =1$.
	Hence, $»\xhh{\widetilde  B^Y_t}{B}(»v,»u_n , e )$ converges to $»\xhh{\widetilde  B^Y_t}{B}(»v,»u, e )$ in $ Y ^*$.
\endproof

\begin{corollary}\label{c4:Psm}
	If {\rm(Ha)} is satisfied, then the mapping $ B _t: V \rightarrow   V ^*$ is pseu\discretionary{-}{}{}domono\discretionary{-}{}{}tone and demicontinuous for all $»t\in T $.
\end{corollary}

{\em Proof}. 
	Due to Corollary~\ref{c4:Ass}, $\xhh{\widetilde  B^X_t}{B}$ and $»\xhh{\widetilde  B^Y_t}{B}$ satisfy the conditions {\rm(A1)} and {\rm(A2)} from~Definition~\ref{d1:Semi}
		as well as {\rm(A3.4)} and $\alpha_A >0$, since 
	\[
		\alpha_A \alpha_B  \enspace \ge\enspace   \alpha_A \alpha_B -\beta_A \beta_B  «= (\alpha_{ b_1 ,p} - C_P \,\beta_{ b_2 ,p})\alpha_{ b_0 ,e}-(\beta_{ b_1 ,e}+ C_P \,\beta_{ b_2 ,e})\beta_{ b_0 ,p} «>0. 
	\]
	Moreover, Lemma~\ref{p4:Ssc} implies {\rm(A3.1)}--{\rm(A3.3)}.
	Then, the assertion follows from Theorem~\ref{t1:SemReduct}, Proposition~\ref{p1:SemiPsm} and Proposition~\ref{p1:Demi}.
\endproof

\begin{proposition}\label{p4:Psm2}
	Under condition {\rm(Ha)}, the operator ${\cal B}:{\cal  V }\times{\cal  V }^*$ is bounded, demicontinuous and pseu\discretionary{-}{}{}domono\discretionary{-}{}{}tone with respect to ${\cal  L }$
		and coercive with respect to $0\in{\cal  V }$.
\end{proposition}

{\em Proof}. 
	Remark~\ref{r4:E} guarantees the injectivity of $ K $.
	Therefore, we identify $ V $ with $ K ( V )$ as in	Remark~\ref{r3:E}
		and prove the assertion by showing that the hypotheses of
		\cite[Prop.~1, p.~440]{papa1} are fulfilled.
	The measurability of $»t\mapsto  B (»t,»u)$ and the growth conditions follow from Lemma~\ref{p4:Growth},
		the pseudomonotonicity and the demicontinuity of $»u\mapsto  B (»t,»u)$ from Corollary~\ref{c4:Psm}.
	It therefore remains to show that there are $C>0$ and $»g\in L^1( T )$ with
		\[
			 \iKlammerB{ B (»t,»u) }{»u}  \iSub{ V } \enspace \ge\enspace   »g(»t) + C\, ||\iBlock{»u}||  _{ V } ^2
		\]
		for all $»t\in T ,»u\in V $.
	Using Lemma~\ref{p4:Growth} and Lemma~\ref{p1:LipMon} we obtain
	\begin{eqnarray*} 
			 \iKlammerB{ S _t »u}{»u}  \iSub{ V }  &  =  &   \iKlammerB{\xwt{ S }_t (»u,»u)-\xwt{ S }_t (0,»u)}{»u-0}  \iSub{ V }  »+   \iKlammerB{\xwt{ S }_t (»0,»u)}{»u}  \iSub{ V }  \\[0.5ex]
							& \>\ge\>   &  \f{\alpha_A \alpha_B -\beta_A \beta_B }{\alpha_B }\>  ||\iBlock{»u}||  _{ V } ^2  »-  ||\iBlock{\xwt{ S }_t (»0,»u)}||  _{ V ^*}  ||\iBlock{»u}||  _{ V }  \\[0.5ex]
							& \>\ge\>   &  \f{\alpha_A \alpha_B -\beta_A \beta_B -\alpha_B \varphi _a }{\alpha_B }\>  ||\iBlock{»u}||  _{ V } ^2  »-  \sqrt{»h(»t)}\, ||\iBlock{»u}||  _{ V }  \\[0.5ex]
							& \>\ge\>   &  \f{\alpha_A \alpha_B -\beta_A \beta_B -\alpha_B \varphi _a }{2\alpha_B }\>  ||\iBlock{»u}||  _{ V } ^2  »- C\,»h(»t).
	\end{eqnarray*}
	This shows the coercivity condition and completes the proof.
\endproof

\begin{theorem}[{\rm Existence of weak solution}]\label{}
	If {\rm(Ha)} is satisfied, then there exists a weak solution $(»u,{\bf u})\in{\cal W}\times{\cal U}$ to problem~{\rm(P)}.
\end{theorem}

{\em Proof}. 
	By Corollary~\ref{c4:Iff}, it suffices to show the existence of a solution $»u\in{\cal W}$ to
	\begin{equation} \label{«cc1}
					({\cal  L }+{\cal A}+{\cal B})\,»u \>\ni\> 0,\enspace \enspace \enspace  ({\cal  K }»u)(0)= K u_0 .
	\end{equation} %
	Condition~{\rm(H1a)} implies that ${\cal A}:{\cal  V }\rightarrow  {\cal  V }^*$ is bounded.
	Moreover, together with $\varphi $ also $ Q $ and ${\cal Q}$ are convex, lower-\discretionary{}{}{}semi\discretionary{-}{}{}con\discretionary{-}{}{}tin\discretionary{-}{}{}u\discretionary{-}{}{}ous and proper.
	Hence, $ A $ is maximal monotone.
	By Proposition~\ref{p4:Psm2}, ${\cal B}:{\cal  V }\rightarrow  {\cal  V }^*$ is bounded, demicontinuous, pseu\discretionary{-}{}{}domono\discretionary{-}{}{}tone with respect to ${\cal  L }$ and coercive
		with respect to $0\in D({\cal A})\cap{\cal W}$.
	Therefore, Theorem~\ref{t2:Exist} yields the existence of a solution $»u\in{\cal W}$ to~\nix(\ref{«cc1}).
\endproof

\subsection*{Acknowledgment}

The author gratefully acknowledges the support of the DFG Research Training Group 1128
	"Analysis, Numerics, and Optimization of Multiphase Problems"
	during the time of his dissertation.
The results of the present paper form parts of his PhD thesis~\cite{doni}.



\bibliography{bibl.bib}{}
\bibliographystyle{plain}

\end{document}